\documentclass[a4paper,12pt]{article}

\newtheorem{theorem}{Theorem}

\newtheorem{corollary}{Corollary}
\newtheorem{algorithm}{Algorithm}

\newtheorem{assumption}{Assumption}
\newtheorem{remark}{Remark}

\usepackage{booktabs}

\usepackage{algorithm}
\usepackage[noend]{algpseudocode}
\makeatletter
\def\BState{\State\hskip-\ALG@thistlm}
\makeatother
\usepackage{algorithmicx}

\usepackage{amsmath}
\usepackage{amssymb}
\usepackage{amsfonts}
\usepackage{graphicx}
\usepackage{ dsfont }
\usepackage{hyperref}

\usepackage[normalem]{ulem}
\usepackage{color}

\usepackage{tikz}
\usepackage{lipsum}

\newcommand{\bs}{\boldsymbol}
\newcommand{\mr}[1]{\mathrm{#1}}

\newcommand{\Mc}{\mathcal{M}}
\newcommand{\Ac}{\mathcal{A}}

\newcommand{\Fc}{\mathcal{F}}
\newcommand{\Uc}{\mathcal{U}}
\newcommand{\Kc}{\mathcal{K}}
\newcommand{\Hc}{\mathcal{H}}

\newcommand{\Rb}{\mathbb{R}}
\newcommand{\Nb}{\mathbb{N}}

\newcommand{\argmin}{\operatornamewithlimits{arg\,min}}

\begin{document}

\title{\bf Linear predictors for nonlinear dynamical systems: Koopman operator meets model predictive control}

\author{Milan Korda$^1$, Igor Mezi{\'c}$^1$}

\footnotetext[1]{Milan Korda and Igor Mezi{\'c} are with the University of California, Santa Barbara,\; {\tt milan.korda@engineering.ucsb.edu, mezic@engineering.ucsb.edu}}

\date{Draft of \today}

\maketitle

\setcounter{footnote}{1}

\begin{abstract}

This paper presents a class of linear predictors for nonlinear controlled dynamical systems. The basic idea is to lift (or embed) the nonlinear dynamics into a higher dimensional space where its evolution is approximately linear. In an uncontrolled setting, this procedure amounts to numerical approximations of the Koopman operator associated to the nonlinear dynamics. In this work, we extend the Koopman operator to controlled dynamical systems and apply the Extended Dynamic Mode Decomposition (EDMD) to compute a finite-dimensional approximation of the operator in such a way that this approximation has the form of a linear controlled dynamical system. In numerical examples, the linear predictors obtained in this way exhibit a performance superior to existing linear predictors such as those based on local linearization or the so called Carleman linearization. Importantly, the procedure to construct these linear predictors is completely data-driven and extremely simple -- it boils down to a nonlinear transformation of the data (the lifting) and a linear least squares problem in the lifted space that can be readily solved for large data sets. These linear predictors can be readily used to design controllers for the nonlinear dynamical system using linear controller design methodologies. We focus in particular on model predictive control (MPC) and show that MPC controllers designed in this way enjoy computational complexity of the underlying optimization problem comparable to that of MPC for a linear dynamical system with the same number of control inputs and the same dimension of the state-space. Importantly, linear inequality constraints on the state and control inputs as well as nonlinear constraints on the state can be imposed in a linear fashion in the proposed MPC scheme. Similarly, cost functions nonlinear in the state variable can be handled in a linear fashion. We treat both the full-state measurement case and the input-output case, as well as systems with disturbances / noise. Numerical examples (including a high-dimensional nonlinear PDE control) demonstrate the approach with the source code available online\footnote{Code download:  \url{https://github.com/MilanKorda/KoopmanMPC/raw/master/KoopmanMPC.zip}}.
\end{abstract}


\begin{flushleft}\small
\end{flushleft}

\section{Introduction}
This paper presents a class of linear predictors for nonlinear controlled dynamical systems. By a predictor, we mean an artificial dynamical system that can predict the future state (or output) of a given nonlinear dynamical system based on the measurement of the current state (or output) and current and future inputs of the system. We focus on predictors possessing a \emph{linear} structure that allows established linear control design methodologies to be used to design controllers for nonlinear dynamical systems.

The key step in obtaining accurate predictions of a nonlinear dynamical system as the output of a linear dynamical system is \emph{lifting} of the state-space to a higher dimensional space, where the evolution of this lifted state is (approximately) linear. For uncontrolled dynamical systems, this idea can be rigorously justified using the Koopman operator theory~\cite{mezic2005spectral,mezic2004comparison}. The Koopman operator is a \emph{linear} operator that governs the evolution of scalar functions (often referred to as observables) along trajectories of a given nonlinear dynamical system. A finite-dimensional approximation of this operator, acting on a given finite-dimensional subspace of all functions, can be viewed as a predictor of the evolution of the values of these functions along the trajectories of the nonlinear dynamical system and hence also as a predictor of the values of the state variables themselves provided that they lie in the subspace of functions the operator is truncated on. In the uncontrolled context, the idea of representing a nonlinear dynamical system by an infinite-dimensional linear operator goes back to the seminal works of Koopman, Carleman and Von Neumann~\cite{koopman1931,carleman1932application,koopman1932}. 
The potential usefulness of such linear representations for prediction and control was suggested in~\cite{mezic2004comparison}.


In this work, we extend the definition of the Koopman operator to controlled dynamical systems by viewing the controlled dynamical system as an uncontrolled one evolving on an extended state-space given by the product of the original state-space and the space of all control sequences. Subsequently, we use a modified version of the Extended Dynamic Mode Decomposition (EDMD)~\cite{edmd} to compute a finite-dimensional approximation of this controlled Koopman operator. In particular, we impose a specific structure on the set of observables appearing in the EDMD such that the resulting approximation of the operator has the form of a linear controlled dynamical system. Importantly, the procedure to construct these linear predictors is completely data-driven (i.e., does not require the knowledge of the underlying dynamics) and extremely simple -- it boils down to a nonlinear transformation of the data (the lifting) and a linear least squares problem in the lifted space that can readily solved for large data sets using linear algebra. On the numerical examples tested, the linear predictors obtained in this way exhibit a predictive performance superior compared to both Carleman linearization as well as local linearization methods. For a related work on extending Koopman operator methods to controlled dynamical systems, see~\cite{brunton2016koopman,brunton_dmd_control,proctor2016generalizing,koopman_cont_extend}. See also~\cite{surana_estim,surana_cdc_16} and~\cite{mauroy_id_cdc_16} for the use of Koopman operator methods for state estimation and nonlinear system identification, respectively.

Finally, we demonstrate in detail the use of these predictors for model predictive control (MPC) design; see the survey~\cite{mayne2000constrained} or the book \cite{grune2011nonlinear} for an overview of MPC. In particular, we show that these predictors can be used to design MPC controllers for \emph{nonlinear} dynamical systems with computational complexity comparable to MPC controllers for \emph{linear} dynamical systems with the same number of control inputs and states. Indeed, the resulting MPC scheme is extremely simple: In each time step of closed-loop operation it involves one evaluation of a family of nonlinear functions (the lifting) to obtain the initial condition of the linear predictor and the solution of a \emph{convex} quadratic program affinely parametrized by this lifted initial condition. Importantly, nonlinear cost functions and constraints can be handled in a linear fashion by including all nonlinear terms appearing in these functions among the lifting functions. Therefore, the proposed scheme can be readily used for predictive control of nonlinear dynamical systems, using the tailored and extremely efficient solvers for linear MPC (in our case qpOASES~\cite{ferreau2014qpoases}), thereby avoiding the troublesome and computationally expensive solution of nonconvex optimization problems encountered in classical nonlinear MPC schemes~\cite{grune2011nonlinear}.

The paper is organized as follows: In Section~\ref{sec:basicIdea} we describe the problem setup and the basic idea behind the use of linear predictors for nonlinear dynamical systems. In Section~\ref{sec:rationale} we derive these linear predictors as finite-dimensional approximations to the Koopman operator extended to nonlinear dynamical systems. In Section~\ref{sec:numericalAlg} we describe a numerical algorithm for obtaining these linear predictors. In Section~\ref{sec:MPC} we describe the use of these predictors for model predictive control. In Section~\ref{sec:ext} we discuss extensions of the approach to input output systems (Section~\ref{sec:io}) and to systems with disturbances / noise (Section~\ref{sec:dist}). In Section~\ref{sec:numEx} we present numerical examples.

\section{Linear predictors -- basic idea}\label{sec:basicIdea}
We consider a discrete-time nonlinear controlled dynamical system
\begin{equation}\label{eq:sys}
x^+ = f(x,u),
\end{equation}
where $x \in \Rb^n$ is the state of the system, $u \in \Uc \subset\Rb^m$ the control input, $x^+$ is the successor state and $f$ the transition mapping. The input-output case is treated in Section~\ref{sec:io}.

The focus of this paper is the prediction of the trajectory of~(\ref{eq:sys}) given an initial condition $x_0$ and the control inputs $\{u_0, u_1,\ldots\}$. In particular, we are looking for simple predictors possessing a linear structure which are suitable for linear control design methodologies such as model predictive control (MPC)~\cite{mayne2000constrained}.

\begin{figure*}[htb]
	\begin{picture}(50,100)
		\put(40,0){\includegraphics[width=140mm]{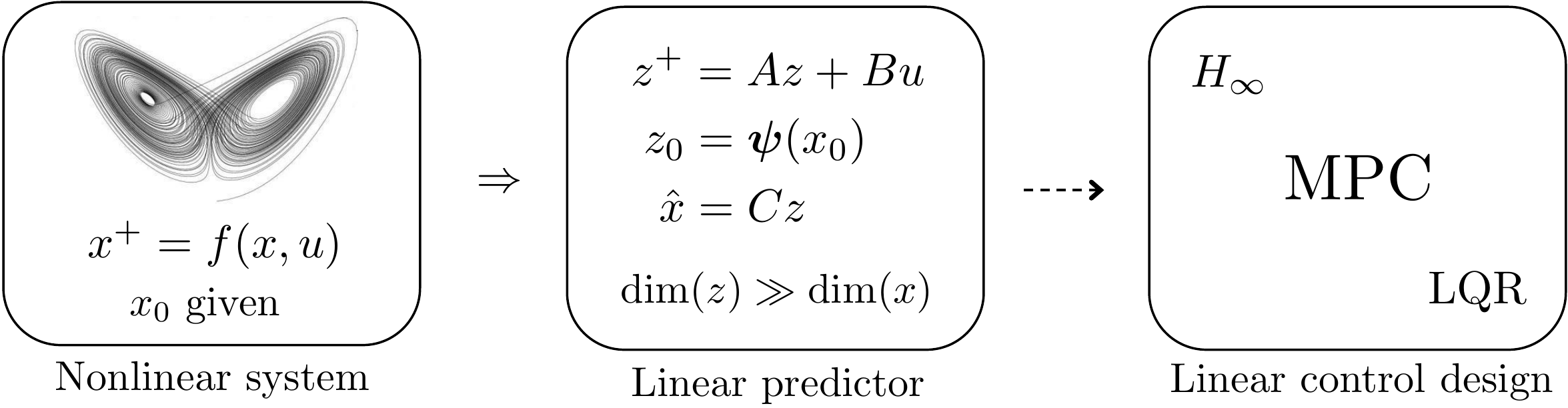}}
	\end{picture}
	\caption{\footnotesize Linear predictor for a nonlinear controlled dynamical system -- $z$ is the lifted state evolving on a higher-dimensional state space, $\hat{x}$ is the prediction of the true state $x$ and $\boldsymbol\psi$ is a nonlinear lifting mapping. This predictor can then be used for control design using linear methods, in our case linear MPC.}
	\label{fig:linPred}
\end{figure*}

The predictors investigated are assumed to be of the form of a controlled \emph{linear} dynamical system
\begin{align}\label{eq:pred_lin}
z^+ & = Az + B u, \\
\hat{x} &= C z \nonumber,
\end{align}
where $z \in \Rb^N$ with (typically) $N \gg n$ and $\hat{x}$ is the prediction of $x$, $B \in \Rb^{N\times m}$ and $C \in \Rb^{n\times N}$. The initial condition of the predictor~(\ref{eq:pred_lin}) is given by
\begin{equation}\label{eq:liftState}
z_0 = \bs\psi(x_0) := \begin{bmatrix} \psi_1(x_0) \\ \vdots \\ \psi_N(x_0) \end{bmatrix} ,
\end{equation}
where $x_0$ is the initial condition of~(\ref{eq:sys}) and $\psi_i : \Rb^n \to \Rb$, $i = 1,\ldots,N$, are user-specified (typically nonlinear) lifting functions. The state $z$ is referred to as the \emph{lifted state} since it evolves on a higher-dimensional, lifted, space\footnote{In general, the term ``lifted state'' may be misleading as, in principle, the same approach can be applied to dynamical systems with states evolving on a (possibly infinite-dimensional) space $\mathcal{M}$, with finitely many observations (or outputs) $\bs h(x) = [h_1(x),\ldots,h_p(x)]^\top$ available at each time instance; the lifting is then applied to these output measurements (and possibly their time-delayed versions), i.e., $\psi_i(x)$ in~(\ref{eq:liftState}) becomes $\psi_i(\bs h(x))$ for some functions $\psi_i : \Rb^p \to \Rb$. In other words, rather than the state itself we lift the output of the dynamical system. For a detailed treatment of the input-output case, see Section~\ref{sec:io}.}. Importantly, the control input $u\in \Uc$ of~(\ref{eq:pred_lin}) remains \emph{unlifted} and hence linear constraints on the control inputs can be imposed in a linear fashion. Notice also that the predicted state $\hat{x}$ is a \emph{linear} function of the lifted state $z$ and hence also linear constraints on the state can be readily imposed. Figure~(\ref{fig:linPred}) depicts this idea.

Predictors of this form lend themselves immediately to linear feedback control design methodologies. Importantly, however, the resulting feedback controller will be nonlinear in the original state $x$ even though it may be (but is not required to be) linear in the lifted state $z$. Indeed, from any feedback controller $\kappa_{
\mr{lift}}: \Rb^N \to \Rb^m$ for~(\ref{eq:pred_lin}), we obtain a feedback controller $\kappa : \Rb^n \to \Rb^m$ for the original system~(\ref{eq:sys}) defined by 
\begin{equation}\label{eq:lift_cont}
\kappa(x) := \kappa_{\mr{lift}}(\bs\psi(x)).
\end{equation}
The idea is that if the true trajectory of $x$ generated by (\ref{eq:sys}) and the predicted trajectory of $\hat{x}$ generated by (\ref{eq:pred_lin}) are close for each admissible input sequence, then the optimal controller for~(\ref{eq:pred_lin}) should be close to the optimal controller for~(\ref{eq:sys}). In Section~\ref{sec:rationale} we will see how the linear predictors~(\ref{eq:pred_lin}) can be derived within the Koopman operator framework extended to controlled dynamical systems.

Note, however, that in general one cannot hope that a trajectory of a linear system~(\ref{eq:pred_lin}) will be an accurate prediction of the trajectory of a nonlinear system for all future times. Nevertheless, if the predictions are accurate on a long enough time interval, these predictors can facilitate the use of linear control systems design methodologies. Especially suited for this purpose is model predictive control (MPC) that uses only \emph{finite-time predictions} to generate the control input.

We will briefly mention in Section~\ref{sec:bilin} how more complex, bi-linear, predictors of the form
\begin{align}\label{eq:pred_bilin}
z^+ & = Az + (Bz) u, \\
\hat{x} &= C z \nonumber
\end{align}
can be obtained within the same framework and argue that predictors of this form can be asymptotically tight (in a well defined sense). Nevertheless, predictors of the from~(\ref{eq:pred_bilin}) are not immediately suited for linear control design and hence in this paper we focus on the linear predictors~(\ref{eq:pred_lin}).


\section{Koopman operator -- rationale behind the approach}\label{sec:rationale}
We start by recalling the Koopman operator approach for the analysis of an \emph{uncontrolled} dynamical system
\begin{equation}\label{sys:uncont}
x^+ = f(x).
\end{equation}
The Koopman operator $\Kc : \Fc \to \Fc$ is defined by
\begin{equation}\label{eq:koopman}
(\Kc \psi)(x) =  \psi(f(x))
\end{equation}
for every $\psi : \Rb^n \to \Rb$ belonging to $\Fc$, which is a space of functions (often referred to as observables) invariant under the action of the Koopman operator. Importantly, the Koopman operator is \emph{linear} (but typically infinite-dimensional) even if the underlying dynamical system is nonlinear. Crucially, this operator fully captures all properties of the underlying dynamical system provided that the space of observables $\Fc$ contains the components of the state $x_i$, i.e, the mappings $x\mapsto x_i$ belong to $\Fc$ for all $i \in \{1,\ldots, n\}$. For a detailed survey on Koopman operator and its applications, see~\cite{applied_koopmanism}.

\subsection{Koopman operator for controlled systems}\label{sec:KoopCont}
There are several ways of generalizing the Koopman operator to controlled systems; see, e.g., \cite{proctor2016generalizing,koopman_cont_extend}. In this paper we present a generalization that is both rigorous and practical. We define the Koopman operator associated to the controlled dynamical system~(\ref{eq:sys}) as the Koopman operator associated to the uncontrolled dynamical system evolving on the extended state-space defined as the product of the original state-space and the space of all control sequences, i.e., in our case $\Rb^n \times \ell(\Uc)$, where $\ell(\Uc)$ is the space of all sequences $(u_i)_{i=0}^\infty$ with $u_i \in \mathcal{U}$. Elements of $\ell({\Uc})$ will be denoted by $\bs u := (u_i)_{i=0}^\infty$. The dynamics of the extended state
\[
\chi = \begin{bmatrix}
x \\ \bs u
\end{bmatrix},
\]
is described by
\begin{equation}\label{eq:sys_ext}
\chi^+ =  F (\chi) := \begin{bmatrix}
f(x,{\bs u}(0)) \\ \mathcal{S} \bs u
\end{bmatrix},
\end{equation}
where $\mathcal{S}$ is the left shift operator, i.e. $(\mathcal{S}\bs{u})(i) = \bs{u}(i+1)$, and $\bs u(i)$ denotes the $i^{\mr{th}}$ element of the sequence $\bs u$.

The Koopman operator $\Kc :  \mathcal{H} \to \mathcal{H}$ associated to~(\ref{eq:sys_ext}) is defined by
\begin{equation}\label{eq:KoopCont}
(\Kc \phi)(\chi) =  \phi( F(\chi))
\end{equation}
for each $\phi:\Rb^n\times \ell(\Uc) \to \Rb$ belonging to some space of observables $\mathcal{H}$.

The Koopman operator~(\ref{eq:KoopCont}) is a linear operator fully describing the non-linear dynamical system~(\ref{eq:sys}) provided that $\mathcal{H}$ contains the components of the non-extended state\footnote{Note that the definition of the Koopman operator implicitly assumes that $\mathcal{H}$ is invariant under the action of $\Kc$ and hence, in the controlled setting, $\Hc$ will typically automatically contain also functions depending on~$\bs u$.} $x_i$, $i=1,\ldots,n$. For example, spectral properties of the operator $\Kc$ should provide information on spectral properties of the nonlinear dynamical system~(\ref{eq:sys}).

\subsection{EDMD for controlled systems}\label{sec:edmd_cont}
In this paper, however, we are not interested in spectral properties but rather in time-domain prediction of trajectories of~(\ref{eq:sys}). To this end, we construct a finite-dimensional approximation to the operator $\Kc$ which will yield a predictor of the form~(\ref{eq:pred_lin}). In order to do so, we adapt the extended dynamic mode decomposition algorithm (EDMD) of~\cite{edmd} to the controlled setting. The EDMD is a data-driven algorithm to construct finite-dimensional approximations to the Koopman operator. The algorithm assumes that a collection of data $(\chi_j, \chi_j^+)$, $j = 1,\ldots,K$ satisfying $\chi_j^+ = F(\chi_j)$ is available and seeks a matrix $\mathcal{A}$ (the transpose of the finite-dimensional approximation of $\Kc$) minimizing
\begin{equation}\label{eq:EDMDobj}
\sum_{j=1}^K \|  \bs \phi(\chi_j^+) - \Ac \bs \phi(\chi_j)   \|^2_2,
\end{equation}
 where
\[
\bs\phi(\chi) = \begin{bmatrix} \phi_1(\chi) & \ldots & \phi_{N_\phi}(\chi) \end{bmatrix}^\top
\]
is a vector of lifting functions (or observables) $\phi_i : \Rb^n \times \ell(\Uc) \to \Rb$, $i\in \{1,\ldots, N_\phi\}$. Note that $\chi = (x,\bs u)$ is in general an infinite-dimensional object and hence the objective~(\ref{eq:EDMDobj}) cannot be evaluated in a finite time unless $\phi_i$'s are chosen in a special way.

\subsubsection{Linear predictors} In order to obtain a linear predictor~(\ref{eq:pred_lin}) and a computable objective function in~(\ref{eq:EDMDobj}) we impose that the functions $\phi_i$ are of the form
\begin{equation}\label{eq:adStruct}
\phi_i(x,\bs u) = \psi_i(x) + \mathcal{L}_i(\bs u),
\end{equation}
where $\psi_i : \Rb^n \to \Rb$ is in general nonlinear but $\mathcal{L}_i: \ell(\Uc) \to \Rb$ is linear. Without loss of generality (by linearity and causality) we can assume that $N_\phi = N + m$ for some $N > 0$ and that the vector of lifting functions $\bs\phi = [\phi_1,\ldots, \phi_{N_{\phi}}]^\top$ is of the form
\begin{equation}\label{eq:basisFun_lin}
\bs\phi(x,\bs u) = \begin{bmatrix}\bs\psi(x) \\ \bs u (0)
\end{bmatrix},
\end{equation}
where $\bs\psi = [\psi_1,\ldots \psi_N]^\top$ and $\bs u(0) \in \Rb^m$ denotes the first component of the sequence $\bs u$. Since we are not interested in predicting future values of the control sequence, we can disregard the last $m$ components of each term $\bs\phi(\chi_j^+) - \Ac \bs \phi(\chi_j) $ in~(\ref{eq:EDMDobj}). Denoting $\bar \Ac$ the first $N$ rows of $\Ac$ and decomposing this matrix such that $\bar\Ac = [A,B]$ with $A \in \Rb^{N\times N}$,  $B \in \Rb^{N\times m}$ and using the notation $\chi_j = (x_j,\bs u_j)$ in~(\ref{eq:EDMDobj}),  leads to the minimization problem
\begin{equation}\label{eq:EDMDobj_mod}
\min_{A,B} \; \sum_{j=1}^K \|  \bs \psi(x_j^+) - A \bs \psi(x_j) - B \bs u_j(0)   \|^2_2.
\end{equation}
Minimizing~(\ref{eq:EDMDobj_mod}) over $A$ and $B$ leads to the predictor of the form~(\ref{eq:pred_lin}) starting from the initial condition
\begin{equation}\label{eq:liftState}
z_0 = \bs\psi(x_0) := \begin{bmatrix} \psi_1(x_0) \\ \vdots \\ \psi_N(x_0) \end{bmatrix}.
\end{equation}
The matrix $C$ is obtained simply as the best projection of $x$ onto the span of $\psi_i$'s in a least squares sense, i.e., as the solution to
\begin{equation}\label{eq:EDMDobj_mod_C}
\min_{C}\;\,\sum_{j=1}^K \|   x_j - C \bs \psi(x_j)  \|^2_2.
\end{equation}
We emphasize that~(\ref{eq:EDMDobj_mod}) and (\ref{eq:EDMDobj_mod_C}) are linear least squares problem that can be readily solved using linear algebra.

\begin{remark}\label{rem:output_proj_triv}
Note that the solution to~(\ref{eq:EDMDobj_mod_C}) is trivial if the set of lifting functions $\{\psi_1,\ldots, \psi_N\}$ contains the state observable, i.e., if, after possible reordering, $\psi_i(x) = x_i$ for all $i\in\{1,\ldots,n\}$. In this case, the solution to~(\ref{eq:EDMDobj_mod_C}) is $C = [I, 0]$, where $I$ is the identity matrix of size $n$.
\end{remark}

 The resulting algorithm for constructing the linear predictor~(\ref{eq:pred_lin}) is concisely summarized in Section~\ref{sec:numericalAlg}.

\subsubsection{Bilinear predictors}\label{sec:bilin}
 Predictors with a more complex structure can be obtained by imposing a structure on the functions $\phi_i$ different than the linear structure~(\ref{eq:adStruct}). In particular, bilinear predictors of the form~(\ref{eq:pred_bilin}) can be obtained by requiring that
\[
\phi_i(x,\bs u) = \psi_i(x) + \xi_i(x)\mathcal{L}_i(\bs u)
\]
for some nonlinear functions $\psi_i$, $\xi_i$ and linear operators~$\mathcal{L}_i$.

A bilinear predictor of the form~(\ref{eq:pred_bilin}) can be tight (in the sense of convergence of predicted trajectories to the true ones as the number of basis functions tends to infinity) under the assumption that the discrete-time mapping~(\ref{eq:sys}) comes from a discretization of a continuous-time system and the discretization interval tends to zero and the underlying continuous-time dynamics is input-affine; see Section IV-C of~\cite{surana_cdc_16} for more details. This bilinearity phenomenon is well known from the classical Carleman linearization in continuous time~\cite{carleman1932application}. In this work, however, we focus on linear predictors since they are immediately amenable to the range of mature linear control design techniques. The use of bilinear predictors for controller design is left for further investigation.

\section{Numerical algorithm -- Finding $A$, $B$, $C$}\label{sec:numericalAlg}

 We assume that a set of data
\begin{equation}\label{eq:data}
 \bf X = \begin{bmatrix} x_1,\;\ldots, x_K \end{bmatrix},\; \bf Y = \begin{bmatrix} y_1,\ldots, y_K \end{bmatrix}, \bf U = \begin{bmatrix} u_1,\ldots, u_K \end{bmatrix}
\end{equation}
satisfying the relation $y_i = f(x_i,u_i)$ is available. Note that we do not assume any temporal ordering of the data. In particular, the data is not required to come from one trajectory of~(\ref{eq:sys}).

Given the data $\bf X$, $\bf Y$, $\bf U$ in~(\ref{eq:data}), the matrices $A \in \Rb^{N\times N}$ and $B \in \Rb^{N\times m}$ in~(\ref{eq:pred_lin}) are obtained as the best linear one-step predictor in the lifted space in a least-squares sense, i.e., they are obtained as the solution to the optimization problem
\begin{equation}\label{eq:ls_AB}
\min_{A,B} \| {\bf Y}_{\mr{lift}} - A {\bf X}_\mr{lift} - B{ \bf U} \|_F,
\end{equation}
where
\begin{equation}\label{eq:liftMat}
\bf X_{\mr{lift}} = \begin{bmatrix} \bs\psi(x_1),\;\ldots, \bs\psi(x_K) \end{bmatrix},\; \bf Y_{\mr{lift}} = \begin{bmatrix} \bs\psi(y_1),\ldots, \bs\psi(y_K) \end{bmatrix},
\end{equation}
 with
\begin{equation}\label{eq:liftState}
\bs\psi(x) := \begin{bmatrix} \psi_1(x) \\ \vdots \\ \psi_N(x) \end{bmatrix}
\end{equation}
being a given basis (or dictionary) of nonlinear functions. The symbol  $\|\cdot\|_{F}$ denotes the Frobenius norm of a matrix. The matrix $C \in \Rb^{n\times N}$ is obtained as the best linear least-squares estimate of $\bf X$ given $\bf X_{\mr{lift}}$, i.e., the solution to
\begin{equation}\label{eq:ls_C}
\min_C \| {\bf X} - C {\bf X}_\mr{lift} \|_F.
\end{equation}

The analytical solution to~(\ref{eq:ls_AB}) is
\begin{equation}\label{eq:ls_anal}
[A, B] = \bf Y_\mr{lift} [\bf X_{\mr{lift}}, \bf U]^\dagger,
\end{equation}
where $^\dagger$ denotes the Moore-Penrose pseudoinverse of a matrix. The analytical solution to~(\ref{eq:ls_C}) is
\[
C = \bf X \bf X_{\mr{lift}}^\dagger.
\]

Notice the close relation of the resulting algorithm to the DMD with control proposed in~\cite{brunton_dmd_control}. There, however, no lifting is applied and the the least squares fit~(\ref{eq:ls_AB}) is carried out on the original data, limiting the predictive power.

\subsection{Practical considerations}
The analytical solution~(\ref{eq:ls_anal}) is not the preferred method of solving the least-squares problem~(\ref{eq:ls_AB}) in practice. In particular, for larger data sets with $K \gg N$ it is beneficial to solve instead the normal equations associated to~(\ref{eq:ls_AB}). The normal equations read
\begin{equation}\label{eq:normEq}
\bf V = \mathcal{M}\bf G
\end{equation}
with variable $\mathcal{M} = [A,B]$ and data 
\[
{\bf G} = \begin{bmatrix}{\bf X_{\mr{lift}}} \\ {\bf U}\end{bmatrix}\begin{bmatrix} {\bf X_{\mr{lift}}} \\ {\bf U} \end{bmatrix}^\top\;,\quad  {\bf V}={\bf Y_{\mr{lift}}}  \begin{bmatrix} {\bf X_{\mr{lift}}} \\ {\bf U}\end{bmatrix}^\top  .
\]
Any solution to~(\ref{eq:normEq}) is a solution to~(\ref{eq:ls_AB}). Importantly, the size of the matrices $\bf G $ and $\bf V$ is $(N+m)\times (N+m)$ respectively $N\times (N+m)$ and hence independent of the number of samples $K$ in the data set~(\ref{eq:data}). The same considerations hold for the least-squares problem~(\ref{eq:ls_C}). Note that, in practice, the lifting functions $\psi_i$ will typically contain the state itself in which case the solution to~(\ref{eq:ls_C}) is just the selection of appropriate indices of $\bf X_{\mr{lift}}$, i.e., after possible reordering, $C = [I, 0]$.

If lifting to a very high dimensional space is required, it may be worth exploring the so called kernel methods known from machine learning which do not require an explicit evaluation of the lifting mapping $\bs\psi$. These methods were successfully applied to the standard EDMD algorithm in~\cite{williams2014kernel}, leading to substantial computational savings.

\section{Model predictive control}\label{sec:MPC}
\looseness-1 In this section we describe how the linear predictor~(\ref{eq:pred_lin}) can be used to design an MPC controller for the nonlinear system~(\ref{eq:sys}) with computational complexity comparable to that of an MPC controller for a \emph{linear} system of the same state-space dimension and number of control inputs. We recall that MPC is a control strategy where the control input at each time step of the closed-loop operation is obtained by solving an optimization problem where a user-specified cost function (e.g., the energy or tracking error) is minimized along a prediction horizon subject to constraints on the control inputs and state variables. Traditionally, linear MPC solves a convex quadratic program, thereby allowing for an extremely fast evaluation of the control input. Nonlinear MPC, on the other hand, solves a difficult non-convex optimization problem, thereby requiring far more computational resources and/or relying on local solutions only; see, e.g., \cite{grune2011nonlinear} for an overview of nonlinear MPC. Just as linear MPC, the lifting-based MPC for nonlinear systems developed here relies on \emph{convex} quadratic programming, thereby avoiding all issues associated with non-convex optimization and allowing for an extremely fast evaluation of the control input. We first describe the proposed MPC controller in its most general form and subsequently, in Section~\ref{sec:NMPCtoKMPC}, describe how a traditional nonlinear MPC problem translates to the proposed one.

The proposed model predictive controller solves at each time instance $k$ of the closed-loop operation the optimization problem
\begin{equation}\label{eq:MPC}
\begin{array}{ll}
\underset{u_i, z_i }{\mbox{minimize}} &J\left((u_i)_{i=0}^{N_p-1}, (z_i)_{i=0}^{N_p}  \right)  \\
\mbox{subject to} & z_{i+1} = Az_i + Bu_i, \quad  \hspace{3mm} i = 0,\ldots, N_p-1 \\
& E_i z_i + F_iu_i \le b_i,  
\hspace{1mm}\quad \quad i = 0,\ldots, N_p-1\\
&E_{N_p}z_{N_p} \le b_{N_p}\\
\mbox{parameter} & z_0 = {\bs \psi}(x_k) ,
\end{array}
\end{equation}
where $N_p$ is the prediction horizon and the \emph{convex} quadratic cost function $J$ is given by
\begin{align*}
J\left((u_i)_{i=0}^{N_p-1}, (z_i)_{i=0}^{N_p}  \right)=z_{N_p}^\top Q_{N_p} z_{N_p} + q_{N_p}^\top z_{N_p} + \sum_{i=0}^{N_p-1} z_i^\top Q_i z_i + u_i^\top R_i u_i + q_i^\top z_i + r_i^\top u_i 
\end{align*}
with $Q_i \in \Rb^{N\times N}$ and $R_i \in \Rb^{m\times m}$  positive semidefinite. The matrices $E_i \in \Rb^{n_c \times N}$ and $F_i \in \Rb^{n_c \times m}$ and the vector $b_i \in \Rb^{n_c}$ define state and input polyhedral constraints. The optimization problem~(\ref{eq:MPC}) is parametrized by the current state of the nonlinear dynamical system~$x_k$. This optimization problems defines a feedback controller
\[
\kappa(x_k) = u_0^\star(x_k),
\]
where $u_0^\star(x_k)$ denotes an optimal solution to problem~(\ref{eq:MPC}) parametrized by the current state~$x_k$.

Several observations are in order:
\begin{enumerate}
\item The optimization problem~(\ref{eq:MPC}) is a \emph{convex} quadratic programming problem (QP).
\item At each time step $k$, the predictions are initialized from the lifted state $\bs\psi(x_k)$.
\item Nonlinear functions of the original state $x$ can be penalized in the cost function and included among the constraints by including these nonlinear functions among the lifting functions $\psi_i$. For example, if one wished for some reason to minimize $\sum_{i=0}^{N_p-1} \cos(\| x_i\|_\infty)$, one could simply set $\psi_1 = \cos(\| x\|_\infty)$  and $q = [1, 0, 0, \ldots, 0]^\top$. See Section~\ref{sec:NMPCtoKMPC} for more details.
\end{enumerate}

\subsection{Eliminating dependence on the lifting dimension}
In this section we show that the computational complexity of solving the MPC problem~(\ref{eq:MPC}) can be rendered independent of the dimension of the lifted state $N$. This is achieved by transforming~(\ref{eq:MPC}) in the so-called \emph{dense form}
\begin{equation}\label{eq:MPC_dense}
\begin{array}{ll}
\underset{U \in \Rb^{mN_p}}{\mbox{minimize}} & U^\top H U^\top  + h^\top U + z_0^\top G U  \\
\mbox{subject to} &  L U + M z_0 \le c \\
\mbox{parameter} &  z_0 = \bs\psi(x_k),   
\end{array}
\end{equation}
for some positive-semidefinite matrix $H \in \Rb^{mN_p\times mN_p}$ and some matrices and vectors $h \in \Rb^{mNp}$, $G\in \Rb^{N\times mNp}$, $L\in \Rb^{n_cN_p\times mNp}$, $M\in \Rb^{n_cN_p\times N}$ and $c\in \Rb^{n_c N_p}$. The optimization is over the vector of predicted control inputs $U = [u_0^\top, u_1^\top,\ldots, u_{N_p-1}^\top]^\top$. This ``dense'' formulation can be readily derived from the ``sparse'' formulation~(\ref{eq:MPC}) by solving explicitly for $z_i$'s and concatenating the point-wise-in-time stage costs and constraints; see the Appendix for explicit expressions for the data matrices of~(\ref{eq:MPC_dense}) in terms of those of~(\ref{eq:MPC}).

Notice that, crucially, the size of the Hessian $H$ or the number of the constraints $n_c$ in the dense formulation~(\ref{eq:MPC_dense}) is \emph{independent} of the size of the lift $N$. Hence, once the nonlinear mapping $z_0 = \bs\psi(x_k)$ is evaluated, the computational cost of solving~(\ref{eq:MPC_dense}) is comparable to solving a standard \emph{linear} MPC on the same prediction horizon, with the same number of control inputs and with the dimension of the state-space equal to $n$ rather than $N$. This comes from the fact that the cost of solving an MPC problem in a dense form is independent of the dimension of the state-space, once the data matrices in~(\ref{eq:MPC_dense}) are formed. Importantly, these matrices are \emph{fixed} and \emph{precomputed offline} before deploying the controller (with the exception of the inexpensive matrix-vector multiplication $z_0^\top G$). This is in contrast with other MPC schemes for nonlinear systems where these matrices have to be re-computed at each time step of the closed-loop operation, thereby greatly increasing the computational cost.

The closed-loop operation of the lifting-based MPC can be summarized by the following algorithm, where $U^\star_{1:m}$ denotes the first $m$ components of $U^\star$:

\begin{algorithm}
\caption{Lifting MPC -- closed-loop operation}\label{alg:liftMPC}
\begin{algorithmic}[1]
\For{$k=0,1,\ldots$}
\State Set $z_0 := \bs\psi(x_k)$
\State		Solve~(\ref{eq:MPC_dense}) to get an optimal solution $U^\star$
\State Set $u_k = U^\star_{1:m}$
\State $x_{k+1} = f(x_k,u_k)$ [ $=$ apply $u_k$ to the system~(\ref{eq:sys})]
\EndFor
\end{algorithmic}
\end{algorithm}

\subsection{Transforming NMPC to Koopman MPC}\label{sec:NMPCtoKMPC}
In this section we describe in detail how a traditional nonlinear MPC problem translates\footnote{By translate, we do not mean to rewrite in an equivalent form. The problem~(\ref{eq:MPC}) is of course only an approximation to~(\ref{eq:MPC_NL}) since the lifted linear predictor is not exact unless the dynamical system is linear.} to the proposed MPC~(\ref{eq:MPC}). We assume a nonlinear MPC problem which at each time step $k$ of the closed-loop operation solves the optimization problem
\begin{equation}\label{eq:MPC_NL}
\begin{array}{ll}
\underset{u_i, \bar x_i }{\mbox{minimize}} & l_{N_p}(x_{N_p})+\sum_{i=0}^{N_p-1} l_i(\bar x_i) +  u_i^\top \bar R_i  u_i + \bar r_i^\top u_i  \\
\mbox{subject to} & \bar x_{i+1} = f(\bar x_i,u_i), \quad\hspace{9mm} i = 0,\ldots, N_p-1 \\
& c_{x_i}(\bar x_i) + c_{u_i}^\top u_i \le 0,  
\hspace{0mm}\quad \quad i = 0,\ldots, N_p-1\\
& c_{x_{N_p}}(\bar x_{N_p}) \le 0\\
& \bar x_0 = x_k,
\end{array}
\end{equation}
where the notation $\bar x $ is used to distinguish the predicted state $\bar x$, used only within the optimization problem~(\ref{eq:MPC_NL}), from the true measured state $x$ of the dynamical system~(\ref{eq:sys}). Notice that the true  nonlinear dynamics $x^+ = f(x,u)$ appears as a constraint of~(\ref{eq:MPC_NL}); in addition, the functions $l_i$ and $c_{x_i}$ can be nonlinear and hence the optimization problem~(\ref{eq:MPC_NL}) is in general nonconvex and extremely hard to solve to global optimality.

In order to translate~(\ref{eq:MPC_NL}) to the proposed form~(\ref{eq:MPC}) we assume that a predictor of the form~(\ref{eq:pred_lin}) with matrices $A$ and $B$ has been constructed as described in Section~\ref{sec:numericalAlg}, using a lifting mapping $\bs\psi = [\psi_1,\ldots,\psi_N]^\top$. The matrices $A$, $B$ appear in the first constraint of~(\ref{eq:MPC}) and the lifting mapping $\bs\psi$ is used for initialization $z_0 = \bs\psi(x_k)$. In order to obtain the remaining data matrices of~(\ref{eq:MPC}) we assume without loss of generality that $\psi_{i}(x) = l_i(x)$, $i = 0,\ldots, N_p$ and $\psi_{N_p+i}(x) = c_{x_i}(x)$, $i = 0,\ldots, N_p$. With this assumption, the remaining data is given by $Q_i = 0$, $R_i = \bar{R}_i$, $ r = \bar r_i$, $q_i = [0_{1\times i}, 1, 0_{1\times N-1-i} ]$, $E_i = [0_{1\times N_p+i}, 1, 0_{1\times N-N_p-1-i} ]$, $F_i = c_{u_i}^\top$, $b_i = 0$, where $0_{i\times j}$ denotes the matrix of zeros of size $i\times j$. Note that this derivation assumed that the constraint functions $c_{x_i}$ and $c_{u_{i}}$ are scalar-valued; for vector valued constraint functions, the approach is analogous, setting the lifting functions $\psi_i$ equal to the individual components of the constraint functions $c_{x_i}$.

We also note that this canonical approach always leads to a linear cost function in~(\ref{eq:MPC}). However, in special cases, when some of the cost functions $l_i(x_i)$ are convex quadratic, one may want to use the freedom of the formulation~(\ref{eq:MPC}) and instead of setting $\psi_i = l_{i}$, use the quadratic terms in the cost function of~(\ref{eq:MPC}), thereby reducing the dimension of the lift. See Section~\ref{ex:feedbackBilin} for an example.

\section{Theoretical analysis}
In this section we discuss theoretical properties of the EDMD algorithm of Section~\ref{sec:edmd_cont}. Full exposition of the theoretical analysis of EDMD is beyond the scope of this paper and therefore we only summarize the authors' results obtained concurrently in~\cite{korda_mezic_edmd}, where the reader is referred to for proofs of the theorems stated in this section and additional results. Note that the results of Theorems~\ref{thm:L2proj} and \ref{thm:convFin} are not new and were, to the best of our knowledge, first rigorously proven in~\cite[Section~3.4]{klus2015numerical} and alluded to already in the original EDMD paper~\cite{edmd}; here we state them in a language suitable for our purposes.

We work in an abstract setting with a dynamical system
\[
\chi^+ = F(\chi)
\]
with $F : \Mc \to \Mc$, where $\Mc$ is a given separable topological space. This encompasses both the finite-dimensional setting with $\Mc = \Rb^n$ and the infinite-dimensional controlled setting of Section~\ref{sec:KoopCont} with $\Mc = \Rb^n \times l(\Uc)$. The Koopman operator $\Kc:\Hc\to\Hc$ on a space of observables $\Hc$ (with $\phi:\Mc\to\Rb$ for all $\phi \in \Hc$) is defined as in~(\ref{eq:KoopCont}) by $\Kc\phi = \phi \circ F$ for all $\phi \in \Hc$. We assume the EDMD algorithm of Section~\ref{sec:edmd_cont}, i.e., we solve the optimization problem
\begin{equation}\label{eq:EDMDobj_theor}
\min_{\Ac \in \Rb^{N\times N}}\sum_{i=1}^K \|  \bs \phi(\chi_i^+) - \Ac \bs \phi(\chi_i)   \|^2_2,
\end{equation}
where
\[
\bs\phi(\chi) = \begin{bmatrix} \phi_1(\chi) , \ldots, \phi_{N}(\chi) \end{bmatrix}^\top
\]
with $ \phi_i \in \Hc$, $i=1,\ldots,N$, being linearly independent basis functions.
Denoting $\Hc_N$ the span of $\phi_1,\ldots,\phi_N$, the finite-dimensional approximation of the Koopman operator $\Kc_{N,K} : \Hc_N\to \Hc_N$ obtained from~(\ref{eq:EDMDobj_theor}) is defined for any $g = c^\top \bs\phi \in \Hc_N$ by
\begin{equation}
\Kc_{N,K}g = c^\top \Ac_{N,K}\bs\phi,
\end{equation}
where $\Ac_{N,K}$ is the optimal solution of~(\ref{eq:EDMDobj_theor}).

\subsection{EDMD as $L_2$ projection}
First we give a characterization of the EDMD algorithm as an $L_2$ projection. Given an arbitrary nonnegative measure $\mu$ on $\Mc$, we define the $L_2(\mu)$ projection\footnote{The Hilbert space $L_2(\mu)$ is the space of all square integrable functions with respect to the measure $\mu$.} of a function $g$ onto $\Hc_N$ as
\begin{align}\nonumber
P_N^\mu g &= \argmin_{h\in \Hc_N} \|h - g\|_{L_2(\mu)}= \argmin_{h\in \Hc_N} \int_{\Mc} (h - g)^2\,d\mu\\ &=\bs\phi^\top\argmin_{c\in \Rb^N} \int_{\Mc} (c^\top\bs\phi - g)^2\,d\mu. \label{eq:projDef}
\end{align}

We have the following characterization of $\Kc_{N,K} $.
\begin{theorem}\label{thm:L2proj}
Let $\hat \mu_K$ denote the empirical measure associated to the points $\chi_1,\ldots, \chi_K$, i.e., $\hat\mu_K = \frac{1}{K}\sum_{i=1}^K \delta_{\chi_i}$, where $\delta_{\chi_i}$ denotes the Dirac measure at $\chi_i$. Then for any $g \in \Hc_N$
\begin{equation}\label{eq:KoopL2proj}
\Kc_{N,K}g = P_{N}^{\hat{\mu}_K} \Kc g = \argmin_{h \in \Hc_N} \| h - \Kc g  \|_{L_2(\hat\mu_K )} ,
\end{equation}
i.e.,
\begin{equation}\label{eq:KoopL2_altern}
\Kc_{N,K} = P_N^{\hat\mu_K} \Kc_{|\Hc_N},
\end{equation}
where $\Kc_{|\Hc_N}:\Hc_N\to\Hc$ is the restriction of the Koopman operator to the subspace $\Hc_N$.
\end{theorem}
In words, the operator $\Kc_{N,K}$ is the $L_2$ projection of the operator $\Kc$ on the span of $\phi_1,\ldots,\phi_N$ with respect to the empirical measure supported on the samples $\chi_1,\ldots,\chi_K$.

\subsection{Convergence of EDMD}
Now we turn into analyzing convergence $\Kc_{N,K}$ to $\Kc$ as the number of samples $K$ and the number of basis functions $N$ tend to infinity. 

First we analyze convergence as $K$ tends to infinity. For this  we assume a probabilistic sampling model. That is, we assume that the space $\Mc$ is endowed with a probability measure $\mu$ and that the samples $\chi_1,\ldots,\chi_K$ are independent identically distributed (iid) samples from the distribution $\mu$ and we assume that $\Hc = L_2(\mu)$. We invoke the following non-restrictive assumption
\begin{assumption}[$\mu$ independence]\label{as:indep}
The basis functions $\phi_1,\ldots,\phi_N$ are such that \[\mu\big(\{x\in \Mc \mid c^\top \bs \phi(x)=0 \}\big) = 0\] for all nonzero $c
\in\Rb^N$.
\end{assumption}
This is a natural assumption ensuring that the measure $\mu$ is not supported on a zero level set of a nontrivial linear combination of the basis functions used.

Finally we define
\begin{equation}
\Kc_N = P_N^\mu \Kc_{|\Hc_N},
\end{equation}
i.e., the $L_2(\mu)$ projection of the Koopman operator $\Kc$ on $\Hc_N$. Then we have the following result:
\begin{theorem}\label{thm:convFin}
If Assumption~\ref{as:indep} holds, then with probability one
\begin{equation}\label{eq:opConvNormFin}
\lim_{K\to\infty}\|\Kc_{N,K} -  \Kc_N \| = 0,
\end{equation}
where $\| \cdot\|$ is any operator norm and
\begin{equation}\label{eq:opConvSpecFin}
\lim_{K\to\infty}\mr{dist}\big(\sigma(\Kc_{N,K}),\sigma(\Kc_N)\big) = 0,
\end{equation}
where $\sigma(\cdot) \subset \mathbb{C}$ denotes the spectrum of an operator and $\mr{dist}(\cdot,\cdot)$ the Hausdorff metric on subsets of $\mathbb{C}$.
\end{theorem}

Theorem~\ref{thm:convFin} says that the EDMD approximations $\Kc_{N,K}$ converge in the operator norm to the $L_2(\mu)$ projection of the Koopman operator onto the span of the basis functions $\phi_1,\ldots,\phi_N$.

Having established convergence of $\Kc_{N,K}$ to $\Kc$ we turn to studying convergence of $\Kc_N$ to $\Kc$. Since the operator $\Kc_{N,K}$ is defined only on $\Hc_N$ we extend it to all of $\Hc$ by precomposing with the projection on $\Hc_N$, i.e., we study convergence of $\Kc_N P_N^\mu:\Hc\to\Hc$ to $\Kc:\Hc\to\Hc$.

 We have the following result:
\begin{theorem}\label{thm:strongConv}
If $(\phi_i)_{i=1}^\infty$ forms an orthonormal basis of $\Hc = L_2(\mu)$ and if $\Kc:\Hc\to\Hc$ is continuous, then the sequence of operators $\Kc_NP_N^\mu=P_N^\mu \Kc P_N^\mu $ converges to $\Kc$ as $N\to \infty$ in the strong operator topology, i.e., for all $g\in\Hc$
\[
\lim_{N\to\infty} \|\Kc_NP_N^\mu g - \Kc g \|_{L_2(\mu)}  = 0.
\]
In particular, if $g\in\Hc_{N_0}$ for some $N_0\in \Nb$, then 
\[
\lim_{N\to\infty} \|\Kc_N g - \Kc g \|_{L_2(\mu)} = 0.
\]
\end{theorem}
Theorem~\ref{thm:strongConv} tells us that the sequence of operators~$\Kc_N P_N^\mu$ converges strongly to $\Kc$. For additional results on spectral convergence of $\Kc_N$ to $\Kc$, weak spectral convergence of $\Kc_{N,N}$ (i.e., with $K=N$) to $\Kc$ and for a method to construct $\Kc_N$ directly, without the need for sampling, see~\cite{korda_mezic_edmd}.

Now we use Theorem~\ref{thm:strongConv} to establish convergence of finite-horizon predictions of a given vector-valued observable $g \in \Hc_{N_0}^{n_g}$. For our purposes, the most pertinent situation is when $g$ is the state observable, i.e., $g(x) = x$ in which case the following result pertains to predictions of the state itself. We let $\mathcal{A}_{N,K}$ denote the solution to~(\ref{eq:EDMDobj_theor}) and set $\mathcal{A}_{N} = \lim_{K\to\infty}\mathcal{A}_{N,K}$. Since $g\in \Hc_{N_0}^{n_g}$, it follows that for every $N\ge N_0$ there exists a matrix $C_N \in \Rb^{n_g\times N}$ such that $g = C_N \bs\phi_N$, where $\bs\phi_N = [\phi_1,\ldots,\phi_N]^\top$. With this notation the following result holds:

\begin{corollary}\label{cor:finitePred_spec}
If the assumptions of Theorem~\ref{thm:strongConv} hold and $g\in \Hc_{N_0}^{n_g}$ for some $N_0\in\Nb$, then for any finite prediction horizon $N_p \in \Nb$
\begin{equation}\label{eq:finitePred_specpec}
\lim_{N\to\infty} \sup_{i\in\{0,\ldots,N_p\}} \int_{\Mc} \big( C_N \Ac_N^i \bs\phi_N   - g\circ F^i \big)^2\,d\mu = 0.
\end{equation}
\end{corollary}
In words, predictions over any \emph{finite} horizon converge in the $L_2(\mu)$ norm. Unfortunately, Corollary~\ref{cor:finitePred_spec} does not immediately apply to the linear predictors designed in~Section~\ref{sec:edmd_cont} as the set of basis functions~(\ref{eq:basisFun_lin}) does not form an orthonormal basis of $\Hc$ because of the special structure of these basis functions which ensures linearity of the resulting predictors. In this setting, one can only prove convergence of $\Kc_N $ to  $P_\infty^
\mu \Kc_{| \Hc_\infty}$, where $P_\infty^
\mu$ is the $L_2(\mu)$ projection onto  $\Hc_\infty := \overline{\{ \phi_i\mid i\in \mathbb{N} \}}$.

\section{Extensions}\label{sec:ext}
In this section, we describe extensions of the proposed approach to input-output dynamical systems and to systems with disturbances / noise.
\subsection{Input-output dynamical systems}\label{sec:io}
In this section we describe how the approach can be generalized to the case when full state measurements are not available, but rather only certain output is measured. To this end, consider the dynamical system
\begin{align}\label{eq:sys_io}
x^+ &= f(x,u), \\
y &= h(x),\nonumber
\end{align}
where $y$ is the measured output and $h:\Rb^n\to \Rb^{n_h}$.

\subsubsection{Dynamics~(\ref{eq:sys_io}) is known}\label{sec:io_Known}
If the dynamics~(\ref{eq:sys_io}) is known, one can construct a predictor of the form~(\ref{eq:pred_lin}) by applying the algorithm of Section~\ref{sec:numericalAlg} to data obtained from simulation of the dynamical system~(\ref{eq:sys_io}). In closed-loop operation, the predictor~(\ref{eq:pred_lin}) is then used in conjunction with a state-estimator for the dynamical system~(\ref{eq:sys_io}). Alternatively, one can design, using linear observer design methodologies, a state estimator directly for the linear predictor~(\ref{eq:pred_lin}). Interestingly, doing so obviates the need to evaluate the lifting mapping $z = \bs \psi(x)$ in closed loop, as the lifted state is directly estimated. This idea is closely related to the use of the Koopman operator for state estimation~\cite{surana_estim,surana_cdc_16}.

\subsubsection{Dynamics~(\ref{eq:sys_io}) is not known}\label{sec:io_dynNotKnown}
If the dynamics~(\ref{eq:sys_io}) is not known and only input-output data is available, one could construct a predictor of the form (\ref{eq:pred_lin}) or (\ref{eq:pred_bilin}) by taking as lifting functions only functions of the output $y$, i.e., $\psi_i(x) = \phi_i(h(x))$. This, however, would be extremely restrictive as this severely restricts the class of lifting functions available. Indeed, if, for example, $h(x) = x_1$, only functions of the first component are available. However, this problem can be circumvented by utilizing the fact that subsequent measurements, $h(x)$, $h(f(x))$, $h(f(f(x)))$, etc., are available and therefore we can define observables depending not only on $h$ but also on repeated composition of $h$ with $f$. In practice this corresponds to having the lifting functions depend not only on the current measured output but also on several previous measured outputs (and inputs, in the controlled setting). The use of time-delayed measurements is classical in system identification theory (see, e.g.,~\cite{ljung1998system}) but has also appeared in the context of Koopman operator approximation (see, e.g.,~\cite{brunton2016chaos,tu2013dynamic}).

Assume therefore that we are given a collection of data
\[
 {\bf\tilde{X}} = [\bs \zeta_1,\ldots, \bs \zeta_K ],\;\;   {\bf \tilde{Y}} = [\bs \zeta_1^+,\ldots, \bs \zeta_K^+ ],\;\; {\bf \tilde{U}} = [u_1,\ldots,u_K]
\]
where
\begin{align} \label{eq:zetaForm}
\bs\zeta_i &= \begin{bmatrix} y_{i,n_d}^\top & \bar u_{i,n_d-1}^\top & y_{i,n_d-1}^\top & \ldots & \bar u_{i,0}^\top &  y_{i,0}^\top
\end{bmatrix}^\top \in \Rb^{(n_d+1)n_h + n_d m}\\\nonumber
 \bs\zeta_i^+ &= \begin{bmatrix} y_{i,n_d+1}^\top &  \bar u_{i,n_d}^\top & y_{i,n_d}^\top & \ldots & \bar u_{i,1}^\top &  y_{i,1}^\top
\end{bmatrix}^\top \in \ \Rb^{(n_d+1)n_h + n_d m}\\\nonumber
 u_i &= \bar u_{i,n_d}
\end{align}
with $(y_{i,j})_{j=0}^{n_d+1}$ being a vector of consecutive output measurements generated by $(\bar u_{i,j})_{j=0}^{n_d}$ consecutive inputs. We note that there does not need to be any temporal relation between $\bs \zeta_i$ and $\bs \zeta_{i+1}$. If, however, $\bs \zeta_i$ and $\bs \zeta_{i+1}$ are in fact successors, then the matrices $\bf\tilde{X}$ and $\bf\tilde{Y}$ have the familiar (quasi)-Hankel structure known from system identification theory.

Computation of a linear predictor then proceeds in the same way as for the full-state measurement: We lift the collected data and look for the best one-step predictor in the lifted space. This leads to the optimization problem
\begin{equation}\label{eq:ls_AB_io}
\min_{A,B} \big\| {\bf \tilde Y}_{\mr{lift}} - A {\bf \tilde X}_\mr{lift} - B{ \bf \tilde U} \big\|_F,
\end{equation}
where
\begin{equation}\label{eq:liftMat}
\bf \tilde X_{\mr{lift}} = \begin{bmatrix} \bs\psi(\bs \zeta_1),\;\ldots, \bs\psi(\bs \zeta_K) \end{bmatrix},\; \bf \tilde Y_{\mr{lift}} = \begin{bmatrix} \bs\psi(\bs \zeta_1^+),\ldots, \bs\psi(\bs \zeta_K^+) \end{bmatrix}
\end{equation}
and
\[
\bs \psi(\bs \zeta) = \begin{bmatrix} \psi_1( \bs \zeta) \\ \vdots \\ \psi_N( \bs \zeta) \end{bmatrix}
\]
is a vector of real or complex-valued, possibly nonlinear, lifting functions. This leads to a linear predictor in the lifted space
\begin{align}\label{eq:pred_lin_io}
z^+ &= Az + Bu \\
 \hat{y} & = C z, \nonumber
\end{align}
where $\hat{y} $ is the prediction of $y$ and $C$ is the solution to
\begin{equation}\label{eq:outPred_outFeed}
\min_{C} \Big\| [y_{1,n_d},\ldots,y_{K,n_d}] - C{\bf \tilde X}_\mr{lift}  \Big\|_F.
\end{equation}
The predictor~(\ref{eq:pred_lin_io}) starts from the initial condition
\[
z_0 = \bs\psi(\bs \zeta_0),
\]
where
\[
\bs\zeta_0 =  \begin{bmatrix} y_{0}^\top & \bar u_{-1}^\top & y_{-1}^\top & \ldots & \bar u_{-n_d}^\top &  y_{-n_d}^\top
\end{bmatrix}^\top
\]
is the vector of $n_d+1$ most recent output  measurements and $n_d$ input measurements. We remark that the solution to~(\ref{eq:outPred_outFeed}) is trivial provided that the outputs and its delays are included among the lifting functions (e.g., if $\psi_j(\bs\zeta) = \bs\zeta(2j)$, $j=0,\ldots,n_d$); see Remark~\ref{rem:output_proj_triv}.

\begin{remark}[Closed-loop operation] The closed-loop operation of the resulting MPC controller follows the steps of Algorithm~\ref{alg:liftMPC}, only the initialization $z_0 = \bs\psi(x_k)$ in line~2 is replaced by $z_0 = \bs\psi(\bs \zeta_k)$, where $\bs \zeta_k = [y_k^\top, u_{k-1}^\top, y_{k-1}^\top, \ldots, u_{k-n_d}^\top, y_{k-n_d}^\top]^\top$.
\end{remark}

\subsection{Disturbance / Noise propagation}\label{sec:dist}
The approach can be readily extended to systems affected by a disturbance or noise of the form
\begin{equation}\label{eq:sys_dist}
x^+ = f(x,u,w),
\end{equation}
where $w$ is the disturbance or process noise. The goal is to construct a predictor for~(\ref{eq:sys_dist}) of the form
\begin{align}\label{eq:pred_lift_dist}
z^+ &= Az + Bu + Dw \\
 \hat{x} & = C z, \nonumber
\end{align}
starting from the initial condition $z_0 = \bs\psi(x_0)$, where $\bs\psi(\cdot)$ is the lifting mapping defined in~(\ref{eq:liftState}). For example, if $w$ is a stochastic disturbance with known distribution, the linear predictor of the form~(\ref{eq:pred_lift_dist}) can be used to approximately compute the distribution of the state $x$ at a future time instance, given a sequence of control inputs up to that time.

In order to obtain the matrices of the predictor~(\ref{eq:pred_lift_dist}), we assume that  we are given data of the form
\begin{subequations}\label{eq:data_w}
\begin{equation}
 \bf X = \begin{bmatrix} x_1,\;\ldots, x_K \end{bmatrix},\; \bf Y = \begin{bmatrix} y_1,\ldots, y_K \end{bmatrix},
\end{equation}
\begin{equation} 
  \bf U = \begin{bmatrix} u_1,\ldots, u_K \end{bmatrix},\quad \bf W = \begin{bmatrix} w_1,\ldots, w_K \end{bmatrix}
  \end{equation}
\end{subequations}

satisfying $y_i = f(x_i,u_i,w_i)$ for all $i = 1,\ldots,K$. As in Section~\ref{sec:numericalAlg}, the matrices are then obtained using least-squares regression:
\begin{equation}\label{eq:ls_ABD}
\min_{A,B,D} \| {\bf Y}_{\mr{lift}} - A {\bf X}_\mr{lift} - B{ \bf U} - D{ \bf W} \|_F, \quad \min_{C} \| {\bf X} - C {\bf X}_\mr{lift}  \|_F.
\end{equation}
If measurements of the disturbance $w$ are not available, then these must best estimated from the available data, either using one of the nonlinear estimation techniques or using the Koopman operator-based estimator proposed in~\cite{surana_estim,surana_cdc_16}. Alternatively, if the mapping $f$ is known (either analytically or in the form of an algorithm) and an algorithm to draw samples from the distribution of $w$ is available, one can obtain data~(\ref{eq:data_w}) by simulation. 

\begin{remark}If full-state measurements are not available, the approach can be readily combined with the approach of Section~\ref{sec:io_Known} or \ref{sec:io_dynNotKnown}.
\end{remark}

\section{Numerical examples}\label{sec:numEx}
In this section we compare the prediction accuracy of the linear lifting-based predictor~(\ref{eq:pred_lin}) with several other predictors and demonstrate the use of the lifting-based MPC proposed in Section~\ref{sec:MPC} for feedback control of a bilinear model of a motor and of the Korteweg--de Vries nonlinear partial differential equation. The source code for the numerical examples is available from \url{https://github.com/MilanKorda/KoopmanMPC/raw/master/KoopmanMPC.zip}.  

\subsection{Prediction comparison}
In order to evaluate the proposed predictor, we compare its prediction quality with that of several commonly used predictors. The system to compare the predictors on is the classical forced Van der Pol oscillator with dynamics given by
\begin{align*}
    \dot{x}_1 &= 2x_2 \\
    \dot{x}_2 &= -0.8x_1 + 2x_2 - 10x_1^2x_2 + u.
\end{align*}
The predictors compared are:
\begin{enumerate}
\item Predictor based on local linearization of the dynamics at the origin,
\item Predictor based on local linearization of the dynamics at a given initial condition $x_0$,
\item Carleman linearization predictor~\cite{carleman1932application},
\item The proposed lifting-based predictor~(\ref{eq:pred_lin}).
\end{enumerate}

In order to obtain the lifting-based predictor, we discretize the dynamics using the Runge-Kutta four method with discretization period $T_s = 0.01\, \mr{s}$ and simulate 200 trajectories over 1000 sampling periods (i.e., 20\,s per trajectory). The control input for each trajectory is a random signal with uniform distribution over the interval $[-1,1]$. The trajectories start from initial conditions generated randomly with uniform distribution on the unit box $[-1,1]^2$. This data collection process results in the matrices $\bf X$ and $\bf Y$ of size $2\times 2\cdot 10^5$ and matrix $\bf U$ of size $1\times 10^5$. The lifting functions $\psi_i$ are chosen to be the state itself (i.e., $\psi_1 =x_1$, $\psi_2 = x_2$) and 100 thin plate spline radial basis functions\footnote{\label{foot:thinplate}Thin plate spline radial basis function with center at $x_0$ is defined by $\psi(x) = \|x-x_0\|^2\log(\|x-x_0\|)$.} with centers selected randomly with uniform distribution on the unit box. The dimension of the lifted state-space is therefore $N = 102$. 

The degree of the Carleman linearization is set to 14, resulting in the size of the Carleman linearization predictor of $120$ (= the number of monomials of degree less than or equal to 14 in two variables). The $B$ matrix for Carleman linearization predictor is set to $[0, 1,0,\ldots,0]^\top$.

\begin{figure*}[h]
\begin{picture}(140,350)
\put(15,180){\includegraphics[width=67mm]{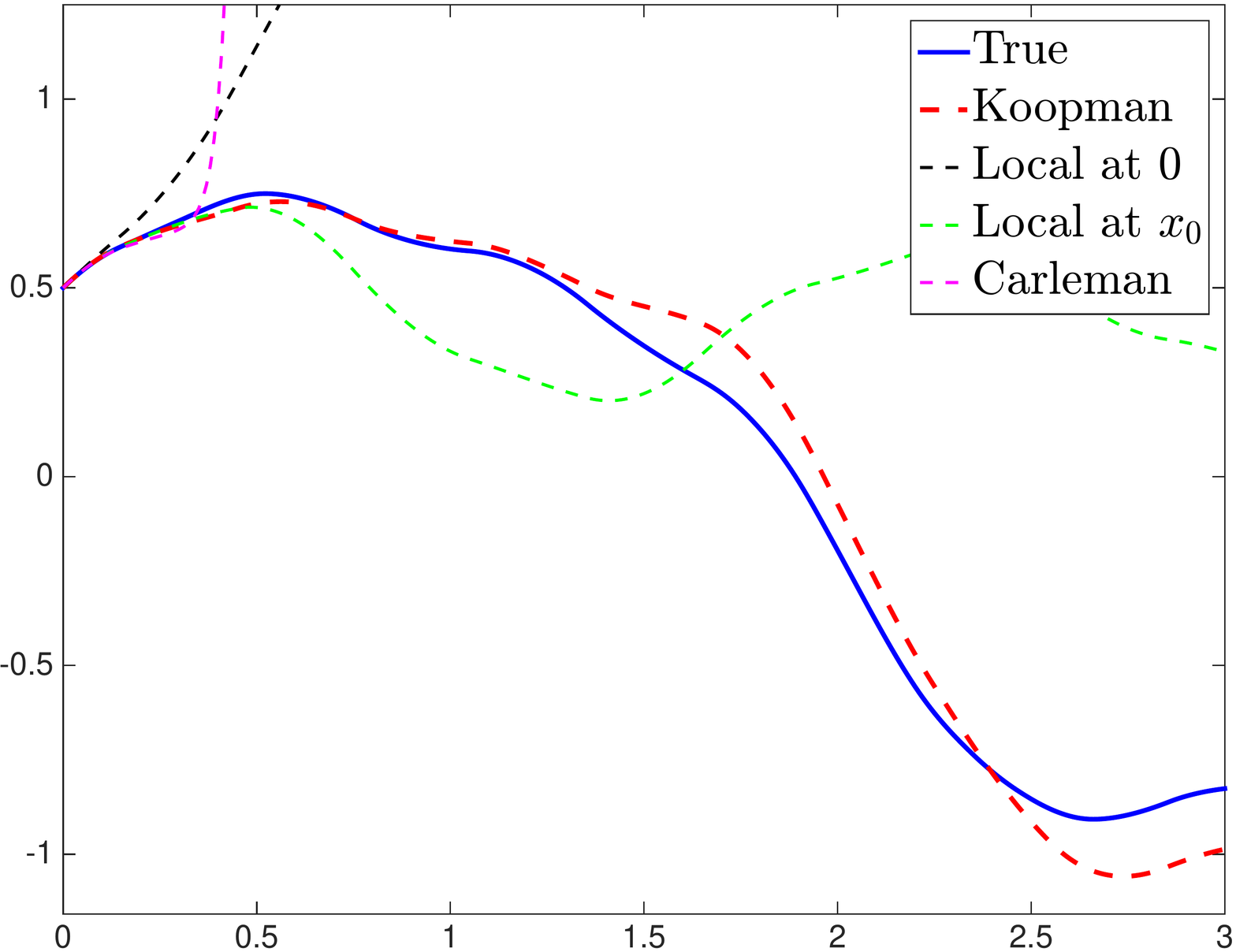}}
\put(245,180){\includegraphics[width=67mm]{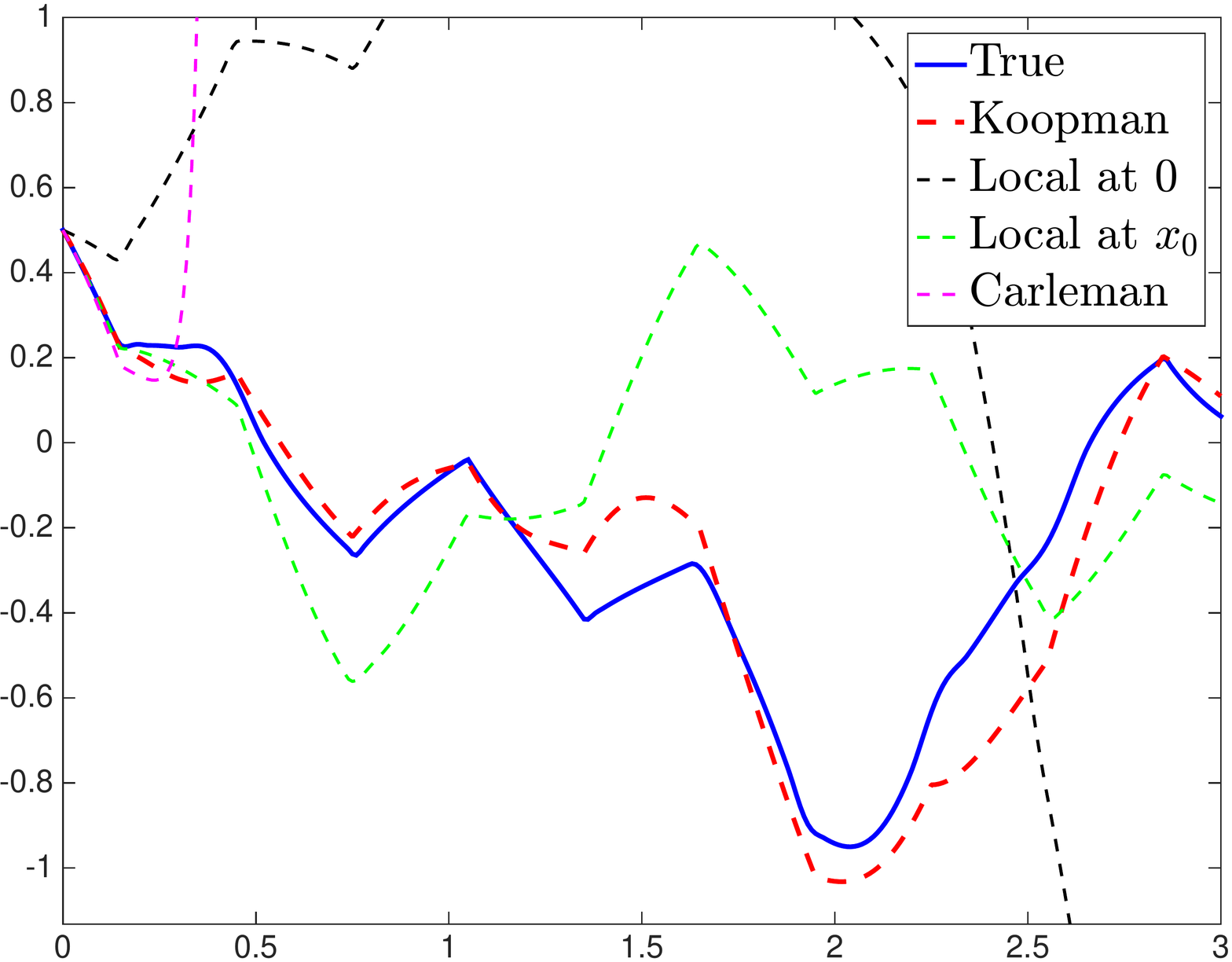}}
\put(15,6){\includegraphics[width=67mm]{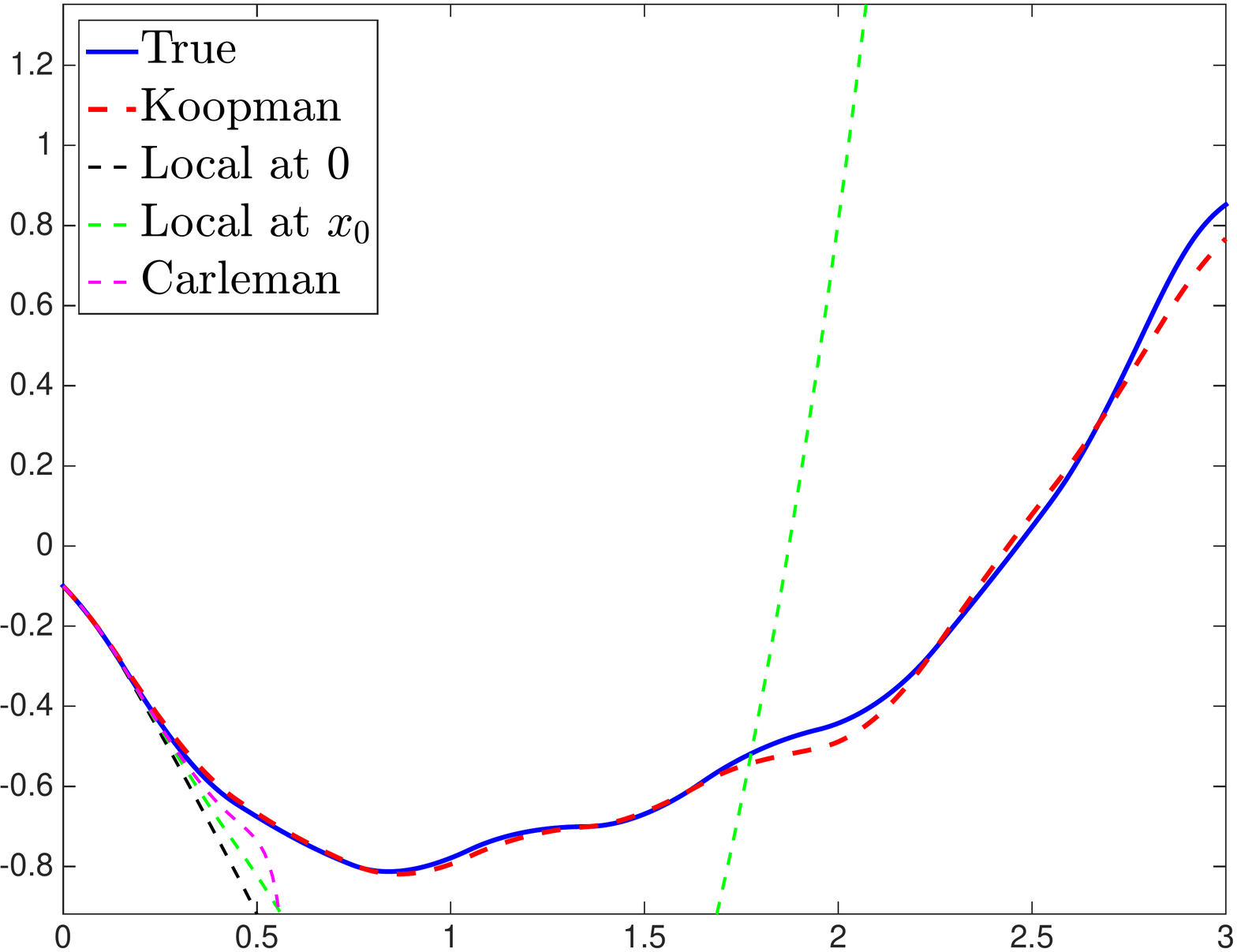}}
\put(245,6){\includegraphics[width=67mm]{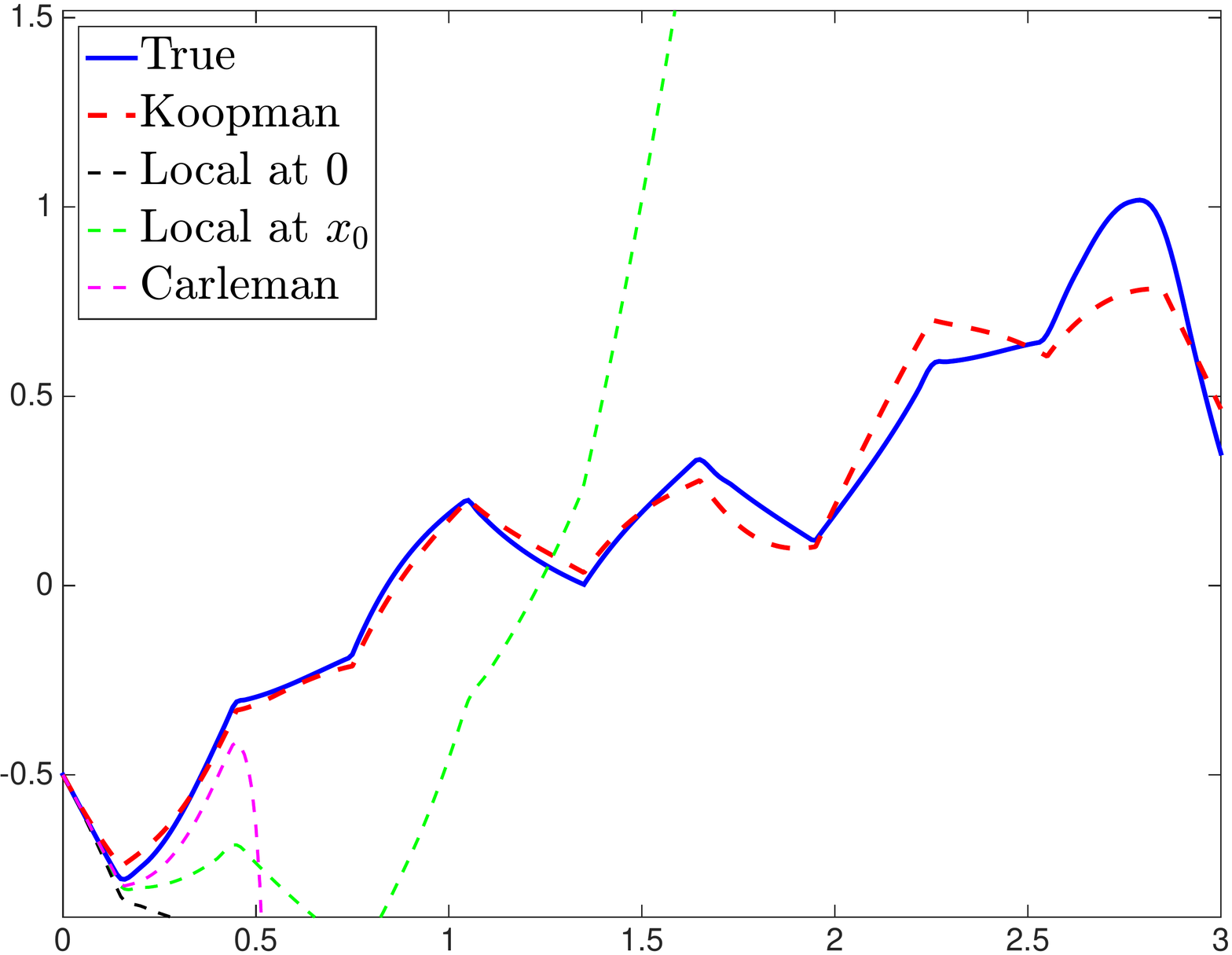}}

\put(100,0){\small $\mr{time}\, [\mr{s}]$}
\put(330,0){\small $\mr{time}\, [\mr{s}]$}
\put(-5,80){ $x_1$}
\put(225,80){ $x_2$}

\put(100,170){\small $\mr{time}\, [\mr{s}]$}
\put(330,170){\small $\mr{time}\, [\mr{s}]$}
\put(-5,250){ $x_1$}
\put(225,250){ $x_2$}

\end{picture}
\caption{Prediction comparison for the forced Van der Pol oscillator. Top: initial condition $x_0 = [0.5,0.5]^\top$. Bottom: initial condition $x_0 = [-0.1,-0.5]^\top$. The forcing $u(t)$ is in both cases a square wave with unit amplitude and period $0.3\,\mr{s}$.}
\label{fig:predictions_vp}
\end{figure*}

Figure~\ref{fig:predictions_vp} compares the predictions starting from two initial conditions $x_0^1 = [0.5,0.5]^\top$, $x_0^2 = [-0.1,-0.5]^\top$ generated by a control signal $u(t)$ being a square wave with unit magnitude and period $0.3\,\mr{s}$. Table~\ref{tab:rmse_vp} reports the relative root mean squared errors (RMSE)
\begin{equation}\label{eq:rmse}
\mr{RMSE} = 100\cdot \frac{\displaystyle\sqrt{\sum_{k}\|x_{\mr{pred}}(kT_s) - x_{\mr{true}}(k T_s)\|_2^2}}{\displaystyle\sqrt{ \sum_{k}\|x_{\mr{true}}(k T_s)\|_2^2}}
\end{equation}
for each predictor averaged over 100 randomly sampled initial conditions with the same square wave forcing. We observe that the lifting-based Koopman predictor is far superior to the remaining predictors. Finally, Table~\ref{tab:rmse_vp_N} reports the prediction accuracy of the lifting predictor in terms of the average RMSE error as a function of the dimension of the lift~$N$; we observe that, as expected, the prediction error decreases with increasing $N$, albeit not monotonously.
 
 \begin{table}[h]
\centering
\caption{\small \rm Prediction comparison --  average RMSE~(\ref{eq:rmse})  over 100 randomly sampled initial conditions: comparison among different predictors. }\label{tab:rmse_vp}\vspace{2mm}
\begin{tabular}{cc}
\toprule
$x_0$ & Average RMSE   \\\midrule
Koopman & 24.4\,$\%$   \\\midrule
Local linearization at $x_0$ & $2.83\cdot 10^3\,\%$ \\ \midrule
 Local linearization at 0  & 912.5\,$\%$  \\\midrule
Carleman  & $5.08\cdot 10^{22}$\,$\%$  \\
\bottomrule
\end{tabular}
\end{table}


\begin{table}[h]
\centering
\caption{\small \rm Prediction comparison -- lifting-based Koopman predictor --  average prediction  RMSE over 100 randomly sampled initial conditions as a function of the dimension of the lift.}\label{tab:rmse_vp_N}\vspace{2mm}
{\small
\begin{tabular}{ccccccc}
\toprule
$N$                            &          5        &        10             &        25                  &           50            &        75               &          100             \\\midrule
  Average RMSE  &   $ 66.5\,\%$   &  $44.9 \,\%$    &           $47.0\,\%$       &       $ 38.7  \,\%$     &   $30.6\,\%$    &      $24.4\,\%$ \\
\bottomrule
\end{tabular}
}
\end{table}

%
%
%

\subsection{Feedback control of a bilinear motor}\label{ex:feedbackBilin}
In this section we apply the proposed approach to the control of a bilinear model of a DC motor~\cite{daniel1998experimental}. The model reads
\begin{align*}
    \dot{x}_1 &= -(R_a/L_a)x_1 - (k_m/L_a)x_2 u + u_a/L_a \\
    \dot{x}_2 &= -(B/J)x_2 + (k_m/J)x_1 u - \tau_l / J, \\
   y  & = x_2
\end{align*}
where $x_1$ is the rotor current, $x_2$ the angular velocity and the control input $u$ is the stator current and the output $y$ is the angular velocity. The parameters are $L_a = 0.314$, $R_a = 12.345$, $k_m = 0.253$, $J = 0.00441$, $B = 0.00732$, $\tau_l = 1.47$, $u_a = 60$. Notice in particular the bilinearity between the state and the control input. The physical constraints on the control input are $u \in [-4,4]$, which we scale to $[-1,1].$

\begin{figure*}[th]
\begin{picture}(140,165)

\put(10,3){\includegraphics[width=75mm]{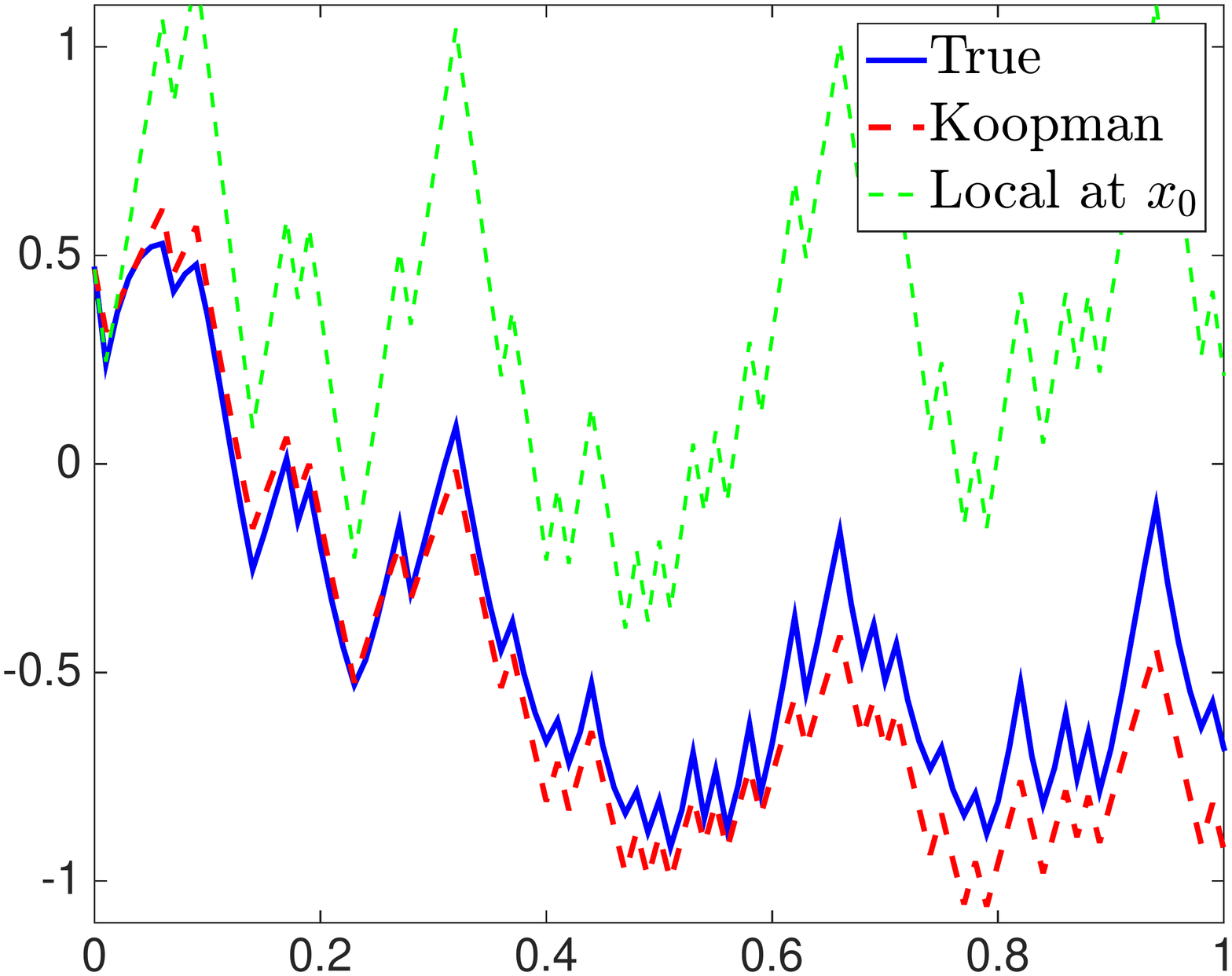}}
\put(230,3){\includegraphics[width=75mm]{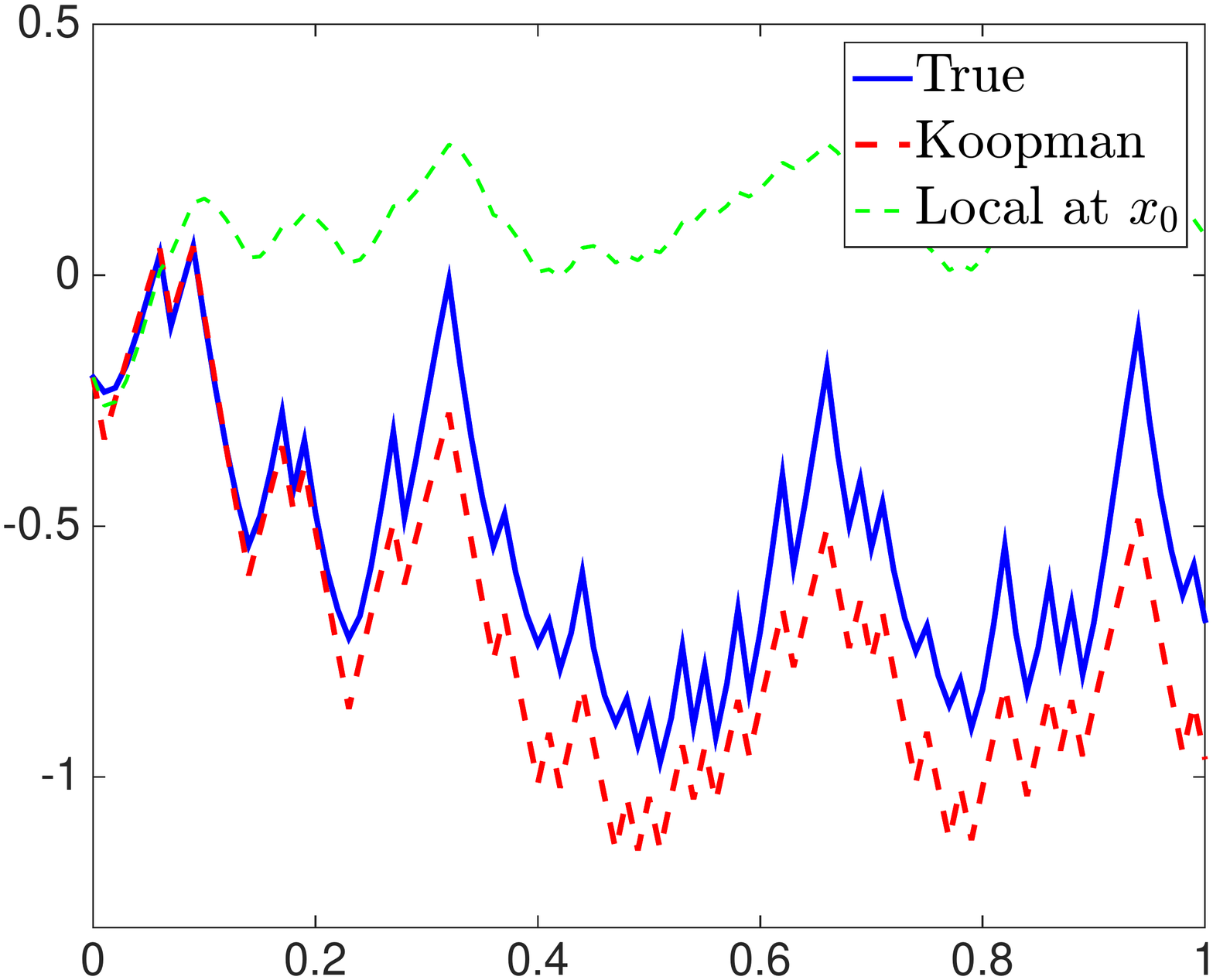}}

\put(105,0){\small $\mr{time} [s]$}
\put(325,0){\small $\mr{time} [s]$}
\put(8,80){ $y$}
\put(228,80){ $y$}
\end{picture}
\caption{\small Feedback control of a bilinear motor -- predictor comparison. Left: initial condition $x_0 =  [ 0.887,0.587]^\top$. Right: initial condition $x_0 =  [-0.404,-0.126]^\top$. }
\label{ref:bilinPredCompar}
\end{figure*}

The goal is to design an MPC controller based on Section~\ref{sec:io_dynNotKnown}, i.e., assuming only input-output data available and no explicit knowledge of the model. In order to obtain the lifting-based predictor~(\ref{eq:pred_lin_io}), we discretize the scaled dynamics using the Runge-Kutta four method with discretization period $T_s = 0.01\, \mr{s}$ and simulate 200 trajectories over 1000 sampling periods (i.e., 20\,s per trajectory). The control input for each trajectory is a random signal uniformly distributed on $[-1,1]$. The trajectories start from initial conditions generated randomly with uniform distribution on the unit box $[-1,1]^2$. We choose the number of delays $n_d = 1$. The lifting functions $\psi_i$ are chosen to be the time-delayed vector $\bs\zeta \in \Rb^3 $, defined in~(\ref{eq:zetaForm}), and 100 thin plate spline radial basis functions (see Footnote~\ref{foot:thinplate}) with centers selected randomly with uniform distribution over $[-1,1]^3$. The dimension of the lifted state-space is therefore $N = 103$. First, in Figure~\ref{ref:bilinPredCompar}, we compare the output predictions for two different, randomly chosen, initial conditions against the predictor based on local linearization at a given initial condition. The prediction accuracy of the proposed predictor is superior, especially for longer prediction times. This is documented further in Table~\ref{tab:predBilinMotor} by the relative root mean-squared errors~(\ref{eq:rmse}) over a one-second prediction horizon averaged over one hundred randomly sampled initial conditions. Both in Figure~\ref{ref:bilinPredCompar} and Table~\ref{tab:predBilinMotor}, the control signal was a pseudo-random binary signal generated anew for each initial condition.


\begin{table}[h]
\centering
\caption{\small \rm Feedback control of a bilinear motor -- prediction RMSE~(\ref{eq:rmse}) for 100 randomly generated initial condtions.}\label{tab:predBilinMotor}\vspace{2mm}
\begin{tabular}{ccc}
\toprule
 & Koopman & Local linearization at $x_0$   \\\midrule
Average RMSE & 32.3\,$\%$ & 135.5\,$\%$ \\
\bottomrule
\end{tabular}
\end{table}

\begin{figure*}[!t]
\begin{picture}(140,300)
\put(20,160){\includegraphics[width=72mm]{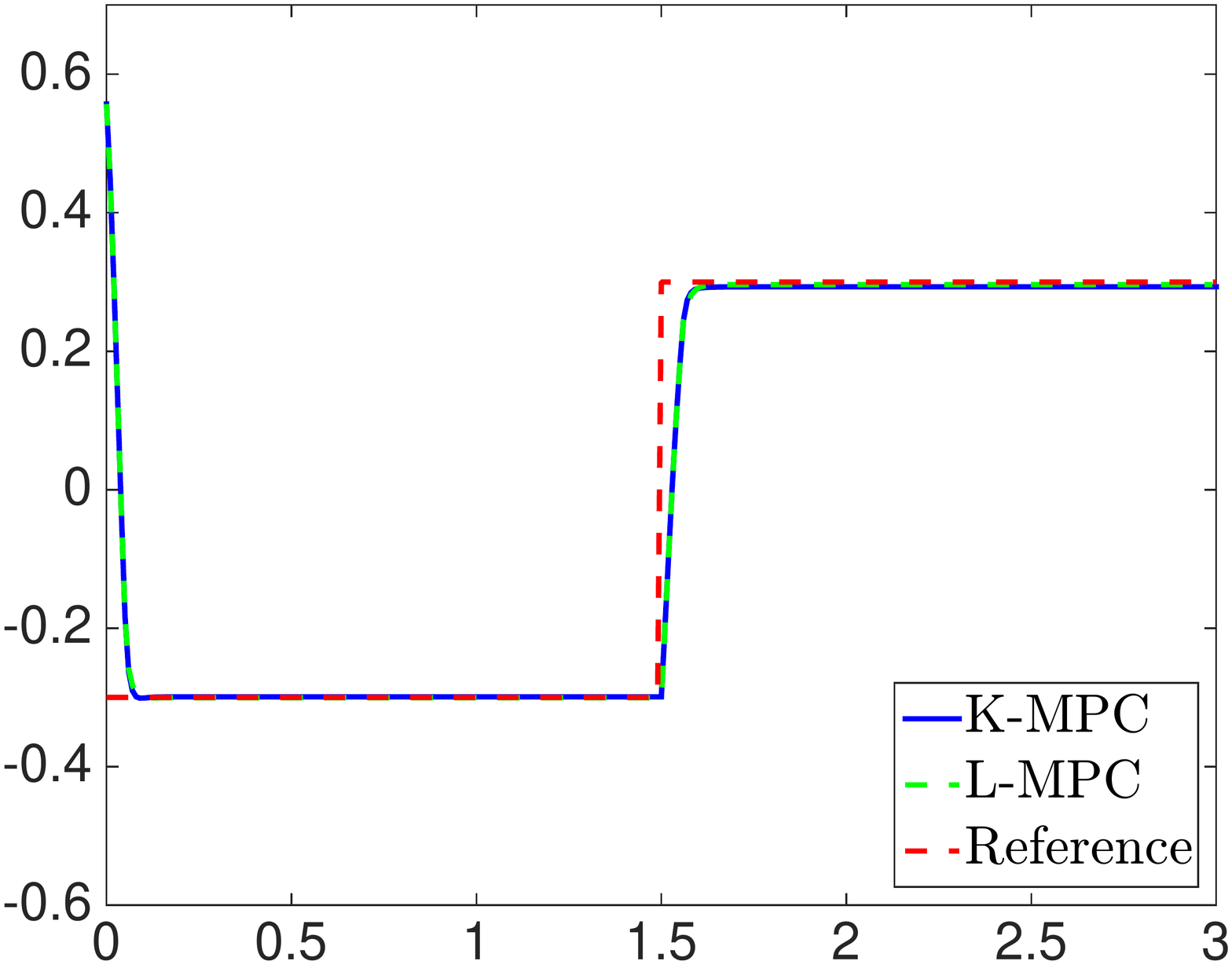}}
\put(220,160){\includegraphics[width=72mm]{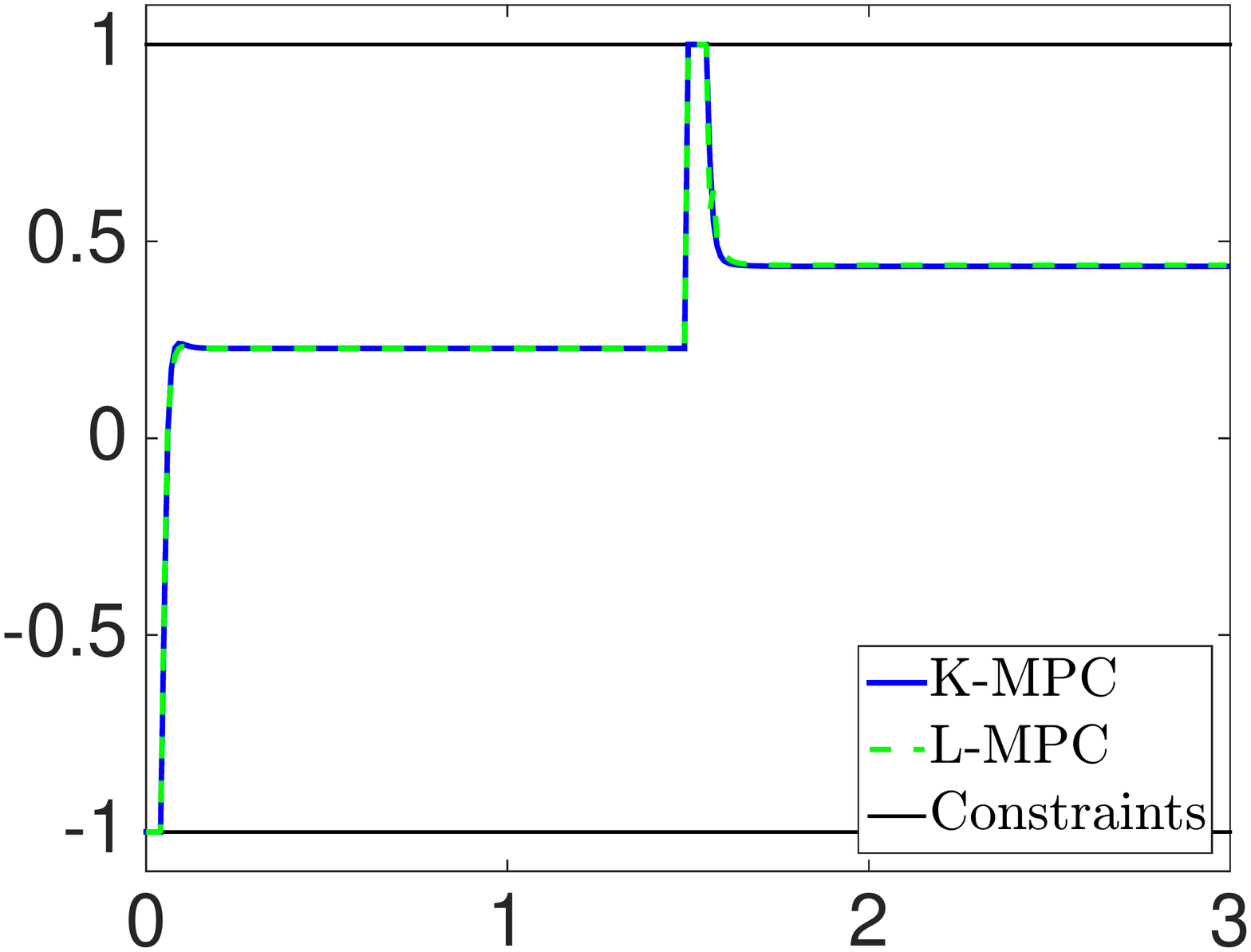}}

\put(20,3){\includegraphics[width=72mm]{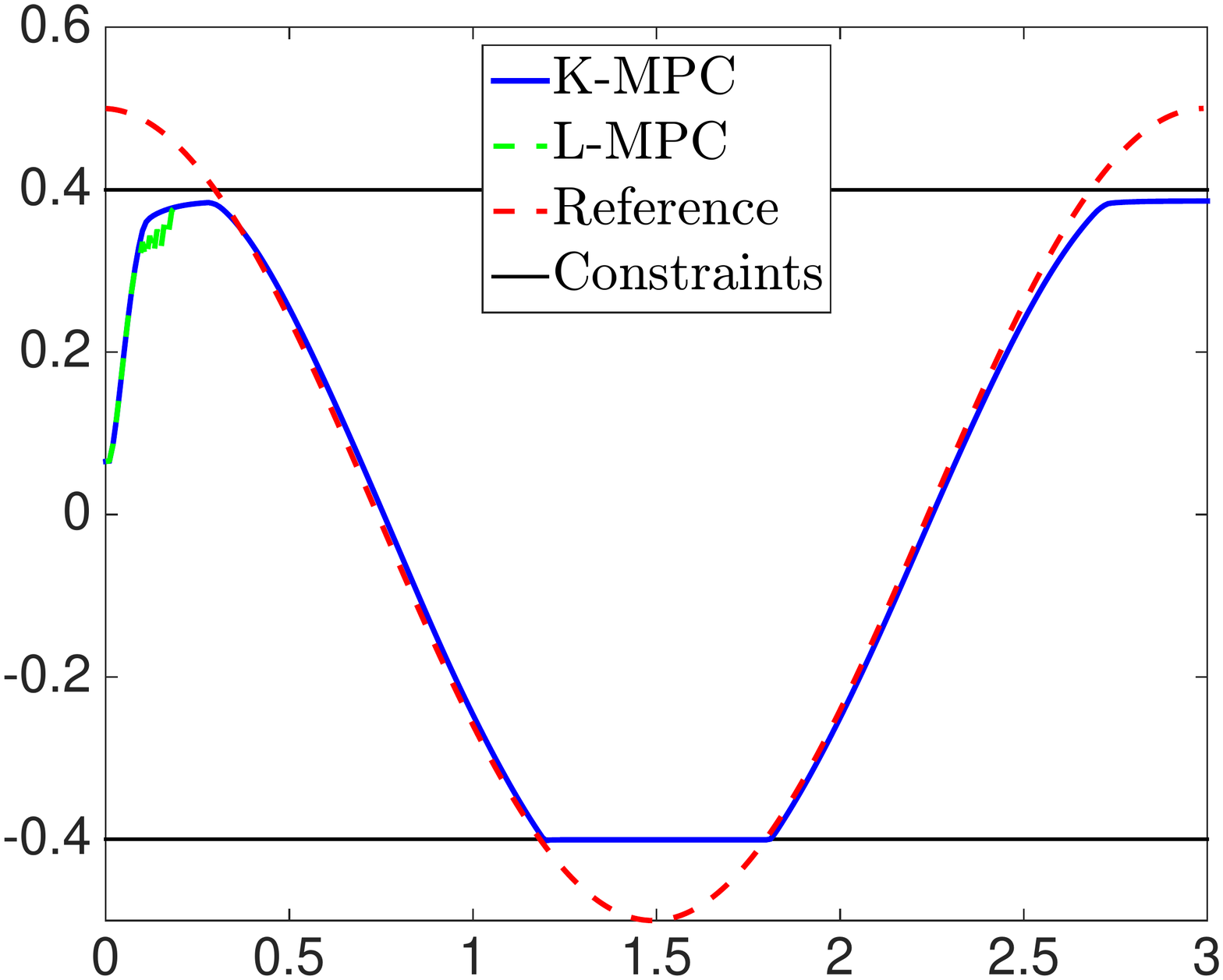}}
\put(220,3){\includegraphics[width=72mm]{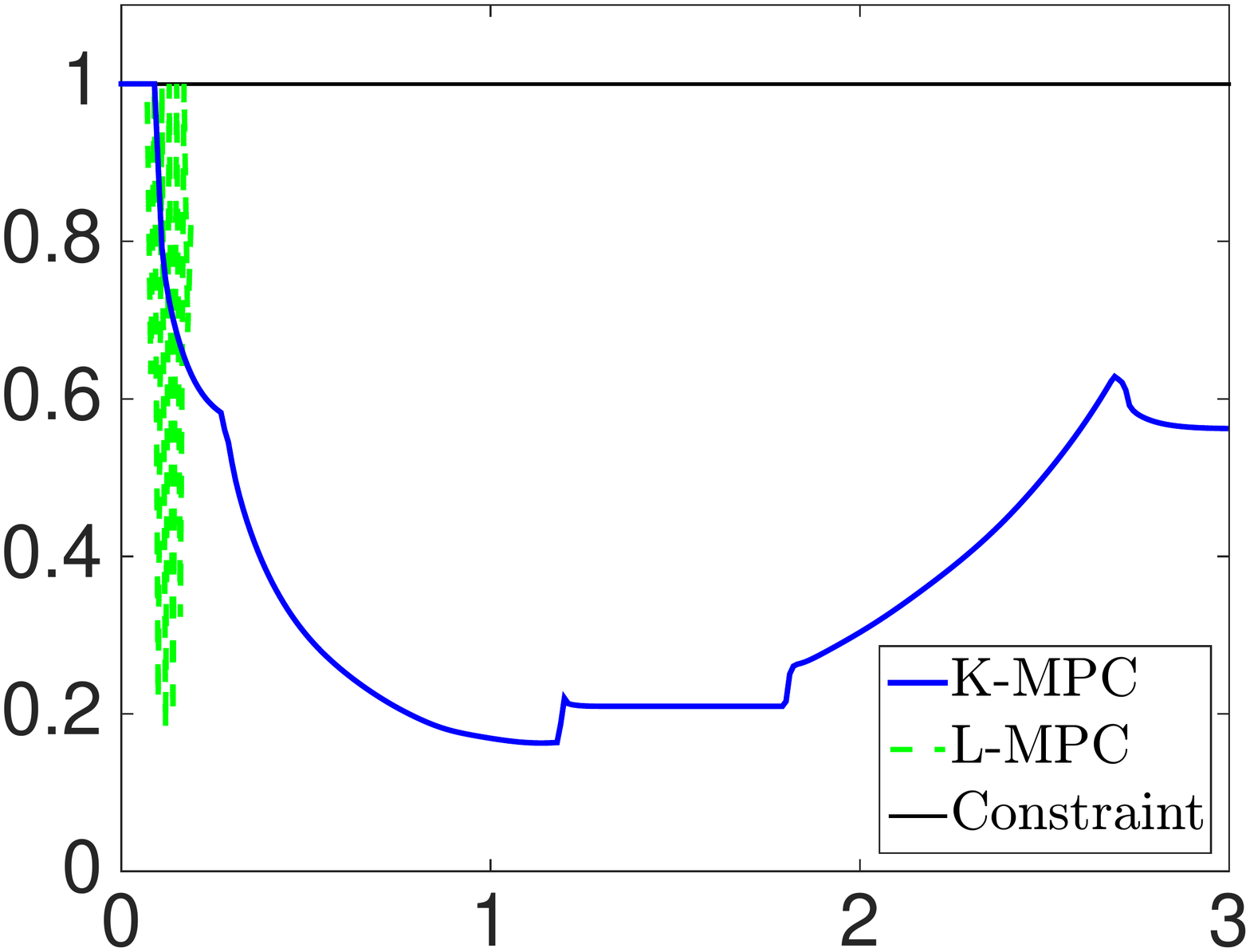}}
\put(271,60){\includegraphics[width=30mm]{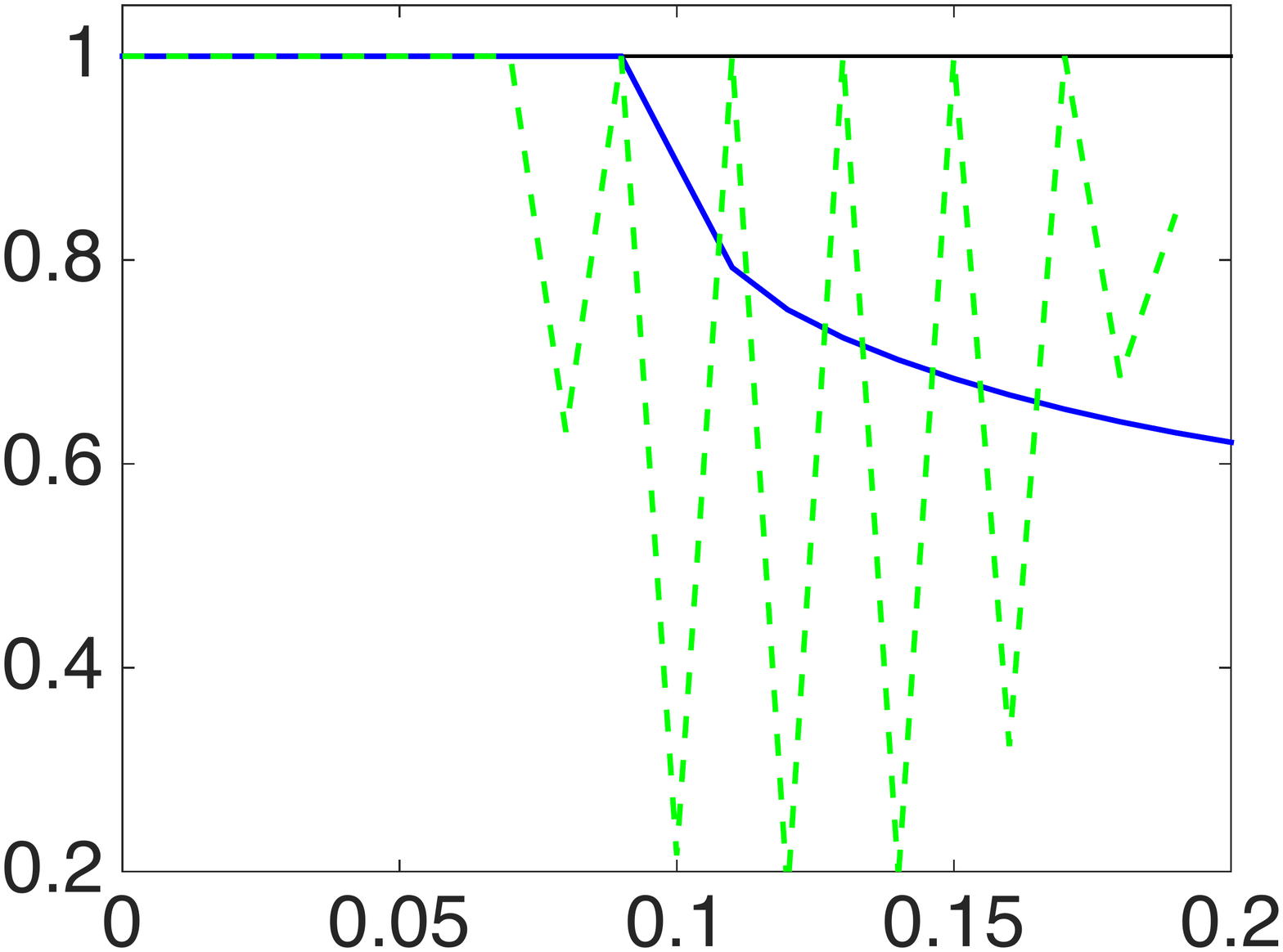}}

\put(68, 65){\vector(-1, 4){13}}
\put(53, 57){\footnotesize{L-MPC}}
\put(50, 47){\footnotesize{infeasible}}

\put(275, 120){\vector(-2, 1){16}}

\put(111,5){\small $\mr{time} [s]$}
\put(311,5){\small $\mr{time} [s]$}
\put(15,80){ $y$}
\put(222,80){ $u$}

\put(110,163){\small $\mr{time} [s]$}
\put(310,163){\small $\mr{time} [s]$}
\put(15,238){ $y$}
\put(222,238){ $u$}

\end{picture}
\caption{\small  Feedback control of a bilinear motor -- reference tracking. Top: piecewise constant reference, $x_0 = [0,0.6]^\top$, no state constraints. Right: time-varying reference, $x_0 = [-0.1,0.1]^\top$, constraints on the output imposed.}
\label{fig:mpc_bilinear_control}
\end{figure*}

The control objective is to track a given angular velocity reference $y_{\mr{r}}$, which translates into the objective function minimized in the MPC problem
\begin{align}
J = &\;\; (C z_{N_p} -y_{\mr{r}})^\top  Q_{N_p} (Cz_{N_p} -y_{\mr{r}})\nonumber \\ & + \sum_{i=0}^{N_p-1} (Cz_{i} -y_{\mr{r}})^\top Q (Cz_{i} -y_{\mr{r}}) + u_i^\top Ru_i\label{obj_track}
\end{align}
with $C = [1,0,\ldots,0]$. This tracking objective function readily translates to the canonical form~(\ref{eq:MPC}) by expanding the quadratic forms and neglecting constant terms. The cost function matrices were chosen as $Q = Q_{N_p} = 1$ and $R = 0.01$. The prediction horizon was set to one second, which results in $N_p = 100$. We compare the Koopman operator-based MPC controller (K-MPC) with an MPC controller based on local linearization (L-MPC) in two scenarios. In the first one we do not impose any constraints on the output and track a piecewise constant reference. In the second one, we impose the constraint $y \in [-0.4,0.4]$ and track a time-varying reference $y_{\mr{r}}(t) = 0.5\cos(2\pi t/3)$, which violates the output constraint for some portion of the simulated period. The simulation results are shown in Figure~\ref{fig:mpc_bilinear_control}. We observe a virtually identical tracking performance in the first case. In the second case, however, the local-linearization controller becomes infeasible and hence cannot complete the entire simulation period\footnote{Infeasibility of the underlying optimization problem is a common problem encountered in predictive control with various heuristic (e.g., soft constraints) or theoretically substantiated (e.g., set invariance) approaches trying to address them. See, e.g., \cite{grune2011nonlinear,MPCbook} for more details.}. This infeasibility occurs due to the inaccurate predictions of the local linearization predictor over longer prediction horizons. The proposed K-MPC controller, on the other hand, does not run infeasible and completes the simulation period without violating the constraints.

Note that, even in the first scenario where the two controllers perform equally, the K-MPC controller has the benefit of being completely \emph{data-driven} and requiring only \emph{output measurements}, whereas the L-MPC controller requires a model (to compute the local linearization) and full state measurements. In addition, the average computation time\footnote{The optimization problems were solved by qpOASES~\cite{ferreau2014qpoases} running on Matlab and 2\,GHz Intel Core i7 with 8\,GB RAM.\label{foot:hw}} required to evaluate the control input of the K-MPC controller was  $6.86\,\mr{ms}$ (including the evaluation of the lifting mapping $\bs \psi(\bs\zeta)$), as opposed to $103\,\mr{ms}$ for the L-MPC controller. This discrepancy is due to the fact that the local linearization and all data defining the underlying optimization problem that depend on it have to be re-computed at every iteration, which is costly on its own and also precludes efficient warm-starting; on the other hand, all data (except for the initial condition) of the underlying optimization problem of K-MPC  are precomputed offline. In both cases, the computation times could be significantly reduced with a more efficient implementation. However, we believe, that the proposed approach would still be  superior in terms of computational speed.

\subsection{Nonlinear PDE control}
\begin{figure*}[!t]
\begin{picture}(140,110)

\put(-30,-30){\includegraphics[width=75mm]{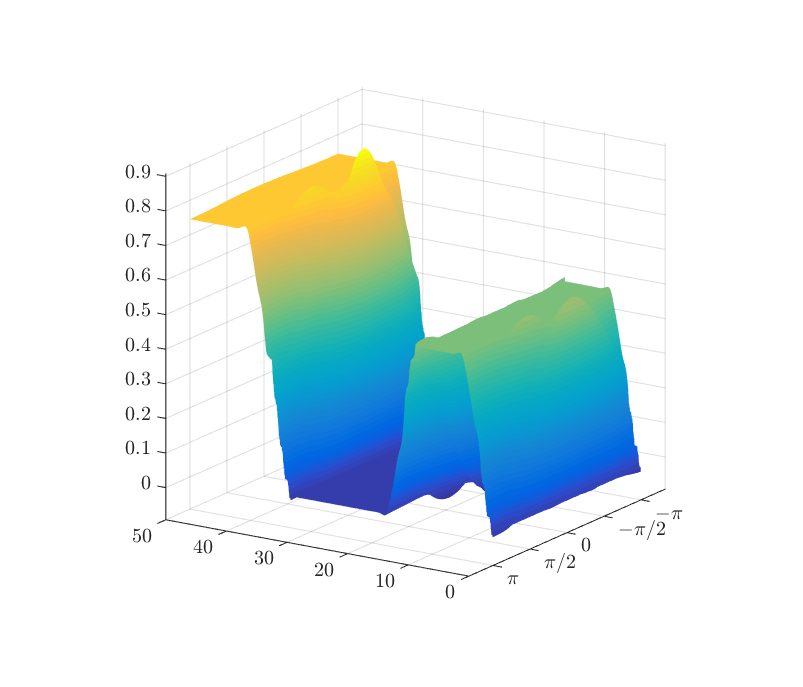}}
\put(170,0){\includegraphics[width=55mm]{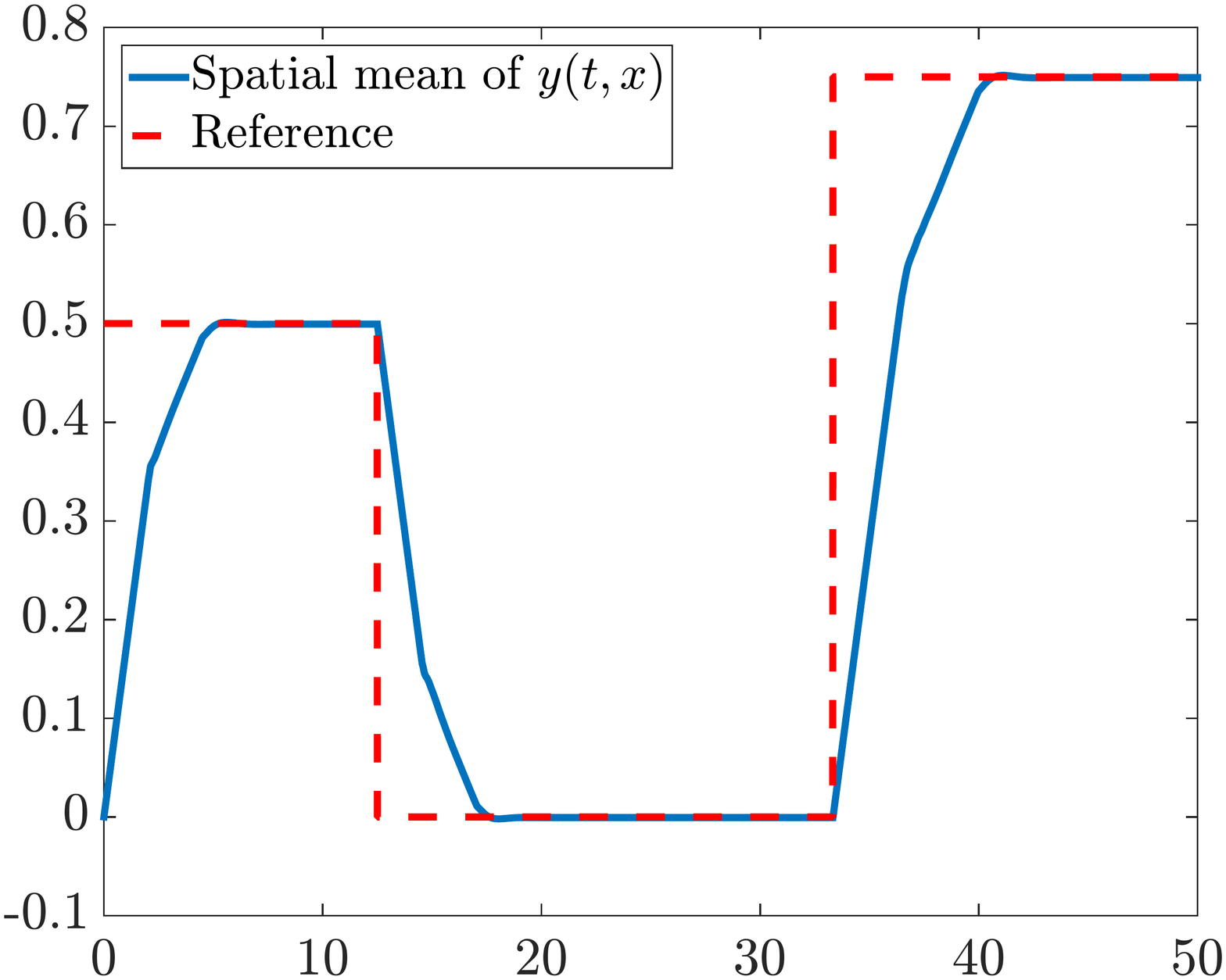}}
\put(330,0){\includegraphics[width=55mm]{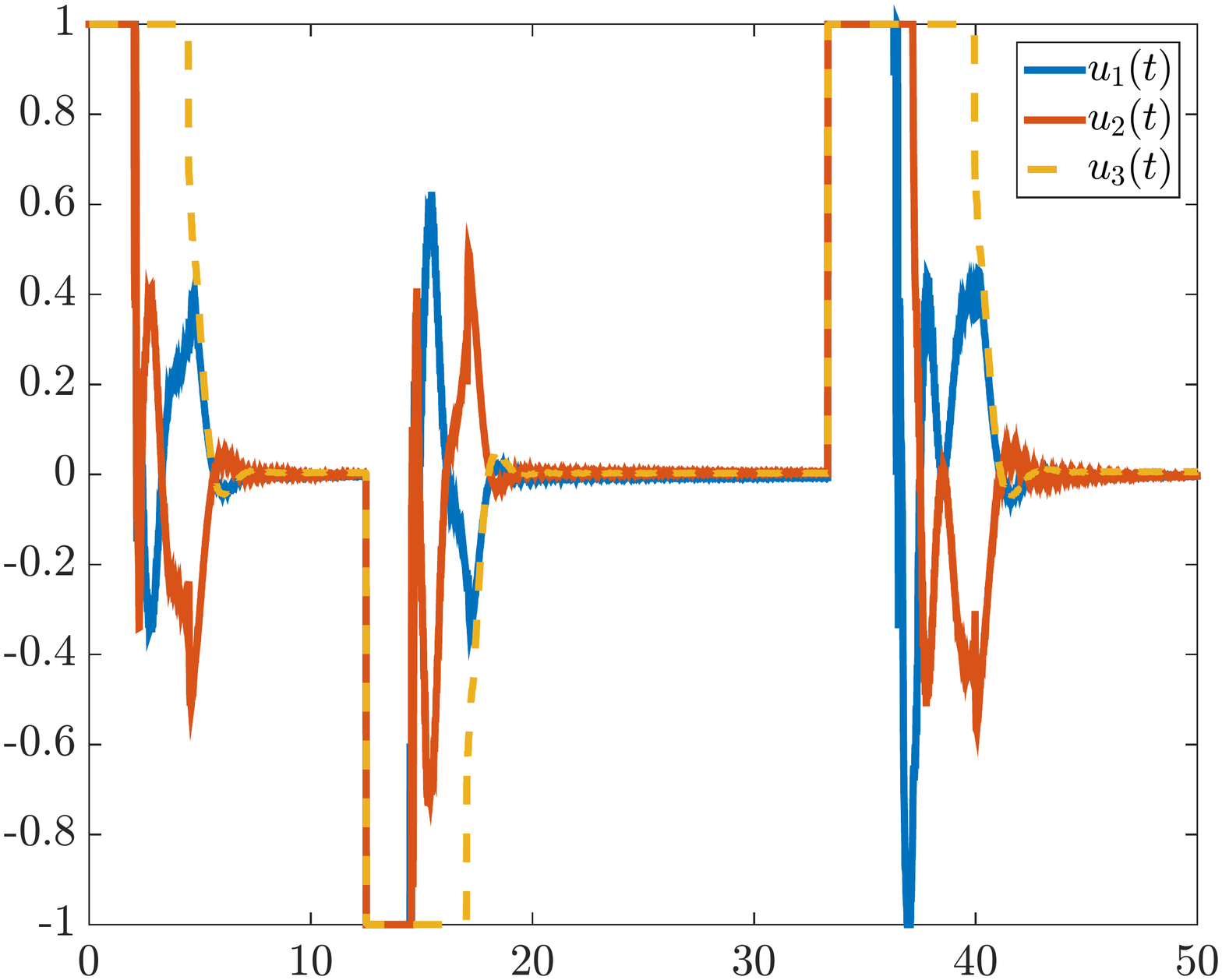}}

\put(128,-3){\footnotesize $x$}
\put(38,-8){\footnotesize $t [\mr{s}]$}
\put(-12,50){\footnotesize \rotatebox{90}{$y(t,x)$}}
\put(245,-5){\footnotesize $t [\mr{s}]$}
\put(402,-5){\footnotesize $t\,[\mr{s}]$}

\end{picture}
\caption{\small Nonlinear PDE control -- Tracking of a time-varying constant-in-space reference profile for the Korteweg--de Vries equation. Left: closed-loop solution. Middle: spatial mean of the solution. Right: control inputs.}
\label{fig:kdv}
\end{figure*}

In order to demonstrate the scalability and versatility of the approach, we use it to control the nonlinear Korteweg--de\,Vries (KdV) equation modelling the propagation of acoustic waves in a plasma or shallow-water waves~\cite{miura1976korteweg}. The equation reads
\[
\frac{\partial\, y(t,x)}{\partial t} + y(t,x)\frac{\partial y(t,x)}{\partial x} + \frac{\partial ^3 y(t,x)}{\partial x^3} = u(t,x),
\]
where $y(t,x)$ is the unknown function and $u(t,x)$ the control input. We consider a periodic boundary condition on the spatial variable $x\in [-\pi,\pi]$. The nonlinear PDE is discretized using the split-stepping method with spatial mesh of $128$ points and time discretization of $\Delta t = 0.01\,\mr{s}$, resulting in a computational state-space of dimension $n=128$. The control input $u$ is considered to be of the form $u(t,x) = \sum_{i=1}^3 u_i(t)v_i(x)$, where the coefficients $u_i(t)$ are to be determined by the controller and $v_i$ are fixed spatial profiles given by $v_i(x) = e^{ -25(x-c_i)^2 }$ with $c_1 = -\pi/2$, $c_2 = 0$, $c_3 = \pi/2$. The control inputs are constrained to $u_i(t) \in [-1,1]$. The lifting-based predictors are constructed from data in the form of 1000 trajectories of length 200 samples. The initial conditions of the trajectories are random convex combinations of three fixed spatial profiles given by $y_0^1 = e^{-(x-\pi/2)^2}$, $y_0^2 = -\sin(x/2)^2$, $y_0^3 = e^{-(x+\pi/2)^2}$; the control inputs $u_i(t)$ are distributed uniformly in $[-1,1]^3$. The lifting mapping $\bs\psi$ is composed of the state itself, the elementwise square of the state, the elementwise product of the state with its periodic shift and the constant function, resulting in the dimension of the lifted state $N = 3\cdot 128 + 1 =385$. The control goal is to track a constant-in-space reference that varies in time in a piecewise constant manner. In order to do so we design the lifting-based Koopman MPC~(\ref{eq:MPC}) with the reference tracking objective~(\ref{obj_track}) with $Q = Q_{N_p} = I$, $R = 0$, $C = [I_{128},0]$ and prediction horizon $N_p = 10$ (i.e., $0.1\,\mr{s}$). The results are depicted in Figure~\ref{fig:kdv}; we observe a fast and accurate tracking of the reference profile. The average computation time to evaluate the control input was $0.28\,\mr{m s}$ (using the dense form~(\ref{eq:MPC_dense}) and the hardware configuration described in Footnote~\ref{foot:hw}), allowing for deployment in real-time applications requiring very fast sampling rates.

\section{Conclusion and outlook}
In this paper, we described a class of linear predictors for nonlinear controlled dynamical systems building on the Koopman operator framework. The underlying idea is to lift the nonlinear dynamics to a higher dimensional space where its evolution is approximately linear. The predictors exhibit superior performance on the numerical examples tested and can be readily used for feedback control design using linear control design methods. In particular, linear model predictive control (MPC) can be readily used to design controllers for the nonlinear dynamical system without resorting to non-linear numerical optimization schemes. Linear inequality constraints on the states and control inputs as well as nonlinear constraints on the states can be imposed in a linear fashion; in addition cost functions nonlinear in the state can be handled in a linear fashion as well. Computational complexity of the underlying optimization problem is comparable to that of an MPC problem for a linear dynamical system of the same size. This is achieved by using the so-called dense form of an MPC problem whose computational complexity depends only on the number of control inputs and is virtually independent of the number of states. Importantly, the entire control design procedure is data-driven, requiring only input-output measurements.

Future work should focus on imposing or proving closed-loop guarantees (e.g., stability or degree of suboptimality) of the controller designed using the presented methodology and on optimal selection of the lifting functions given some prior information on the dynamical system at hand.

\section{Acknowledgments}
The first author would like to thank Colin N. Jones for initial discussions on the topic and for comments on the manuscript, as well as to the three anonymous referees for their constructive comments that helped improve the manuscript. The authors would also like to thank P\' eter Koltai for bringing to our attantion reference~\cite{klus2015numerical}.

 This research was supported in part by the ARO-MURI grant W911NF-14-1-0359 and the DARPA grant HR0011-16-C-0116. The research of M. Korda was supported by the Swiss National Science Foundation under grant P2ELP2\_165166.

\section*{Appendix}\label{ap:matrices}
This appendix expresses explicitly the matrices in the ``dense-form'' MPC problem~(\ref{eq:MPC_dense}) as a function of the data defining the ``sparse-form'' MPC problem~(\ref{eq:MPC}).  The matrices are given by
\[
H = {\mathbf{R}} +{ \mathbf{B}}^\top\mathbf Q {\mathbf{B}},\quad h =  \mathbf B^\top {\mathbf{q}} + \mathbf {r},\quad G = 2\mathbf A^\top\bf Q \mathbf {B},
\]
\[
L = \mathbf{F}+\mathbf{E}\mathbf{B},\quad M = \mathbf{E}\mathbf{A}, \quad c = [b_0^\top,\ldots,b_{N_p}^\top]^\top,
\]
where
{\normalsize
\[
\mathbf{A} = 
\begin{bmatrix}
I \\ A \\ A^2\\ \vdots\\ A^{N_p}
\end{bmatrix}, \,
\mathbf{B} = \begin{bmatrix}0&0&\ldots&0\\
B &0&\ldots&0\\
AB & B & \ldots  & 0 \\
 \vdots &\ddots&\ddots\\
 A^{N_p-1}B & \ldots & AB&B \end{bmatrix},\, \mathbf{F}=\begin{bmatrix}
 F_0& 0 &\ldots &0\\
 0 & F_1 & \ldots & 0\\
 \vdots & & \ddots & \vdots \\
 0 &  0 & \ldots & F_{N_p-1}\\
 0 & 0 & \ldots & 0
 \end{bmatrix}
\]
\[
\mathbf{Q} = \mathrm{diag}(Q_0,\ldots,Q_{N_p}),\;\; \mathbf{R} = \mathrm{diag}(R_0,\ldots,R_{N_p-1}),
\]
\[
\mathbf{E} = \mathrm{diag}(E_0,\ldots,E_{N_p}),\,\mathbf{q} = [q_0^\top,\ldots, q_{N_p}^\top]^\top,\;\; \mathbf{r} = [r_0^\top,\ldots,r_{N_p-1}^\top]^\top,
\]
\normalsize with $\mr{diag}(\cdot,\ldots,\cdot)$ denoting a block-diagonal matrix composed of the arguments.
}

\bibliographystyle{abbrv}
\bibliography{./References}

\end{document}